\documentclass[MSNbibl,number,citesort,dvips]{arxbj}
\usepackage{upgreek}
\usepackage{mathrsfs}
\usepackage{graphicx}

\aid{0}
\volume{18}
\issue{3}
\pubyear{2012}
\firstpage{1002}
\lastpage{1030}
\doi{10.3150/11-BEJ356}

\makeatletter
\DeclareMathAlphabet{\mathbbl}  {U}{mt2hrb}{m}{n}
\SetMathAlphabet{\mathbbl}{bold}{U}{mt2hrb}{b}{n}

\newcommand{\eqref}[1]{(\ref{#1})}

\newtheorem{theorem}{Theorem}
\newtheorem{rmk}{Remark}
\newremark{exe}{Example}

\newcommand{\R}{\mathbb{R}}

\newcommand{\E}{\mathbb{E}}

\newcommand{\calB}{\mathcal{B}}
\newcommand{\e}{\mathrm{e}}
\renewcommand{\P}{\mathbb{P}}
\newcommand{\rt}{\rightarrow}
\newcommand{\X}{\mathbb{X}}
\makeatother

\begin{document}
\begin{frontmatter}

\title{A class of measure-valued Markov chains and Bayesian nonparametrics}
\runtitle{Measure-valued Markov chains}

\begin{aug}
\author[1]{\fnms{Stefano}  \snm{Favaro}\thanksref{1}\ead[label=e1]{favaro@econ.unito.it}},
\author[2]{\fnms{Alessandra} \snm{Guglielmi}\corref{}\thanksref{2}\ead[label=e2]{alessandra.guglielmi@polimi.it}}
\and
\author[3]{\fnms{Stephen G.}~\snm{Walker}\thanksref{3}\ead[label=e3]{S.G.Walker@kent.ac.uk}}
\runauthor{S. Favaro, A. Guglielmi and S.G. Walker}
\address[1]{Universit\`a di Torino and Collegio Carlo Alberto, Dipartimento di Statistica e
Matematica Applicata, Corso Unione Sovietica 218/bis, 10134 Torino, Italy. \printead{e1}}
\address[2]{Politecnico di Milano, Dipartimento di Matematica, P.zza Leonardo da Vinci 32,
20133~Milano, Italy. \printead{e2}}
\address[3]{University of Kent, Institute of Mathematics, Statistics and Actuarial
Science,  Canterbury CT27NZ, UK. \printead{e3}}
\end{aug}

\received{\smonth{2} \syear{2010}}
\revised{\smonth{10} \syear{2010}}

\begin{abstract}
Measure-valued Markov chains have raised interest in Bayesian nonparametrics since the
seminal
paper by (\textit{Math. Proc. Cambridge Philos. Soc.} \textbf{105} (1989) 579--585) 
where a Markov chain having the law of the Dirichlet
process
as unique invariant measure has been introduced. In the present paper, we propose and
investigate a new class of measure-valued Markov chains defined via exchangeable
sequences of random variables.
Asymptotic properties for this new class are derived and applications related to
Bayesian nonparametric mixture modeling, and to a generalization of the Markov chain
proposed by (\textit{Math. Proc. Cambridge Philos. Soc.} \textbf{105} (1989) 579--585), 
are discussed. These results and their
applications highlight once again the interplay between Bayesian nonparametrics and the
theory of measure-valued Markov chains.
\end{abstract}

\begin{keyword}
\kwd{Bayesian nonparametrics}
\kwd{Dirichlet process}
\kwd{exchangeable sequences}
\kwd{linear functionals of Dirichlet processes}
\kwd{measure-valued Markov chains}
\kwd{mixture modeling}
\kwd{P\'olya urn scheme}
\kwd{random probability measures}
\end{keyword}

\end{frontmatter}

\section{Introduction}\label{secintr}

Measure-valued Markov chains, or more generally measure-valued Markov
processes, arise naturally in modeling the composition of evolving populations
and play an important role in a variety of research areas such as population
genetics and bioinformatics (see, e.g., \cite{Eth93,Ethe00,Pit06,Daw93}),
Bayesian nonparametrics \cite{Wal08,Rug08},
combinatorics \cite{Pit06}
and statistical physics \cite{Pit06,Daw93,Delm04}.
In particular, in Bayesian nonparametrics there has been interest in measure-valued Markov chains
since the seminal
paper by \cite{Fei89}, 
where the law of the Dirichlet process has
been characterized as the unique invariant measure of a certain measure-valued
Markov chain.

 In order to introduce the result by \cite{Fei89}, 
let us
consider a Polish space $\mathbb{X}$ endowed with the Borel $\sigma$-field
$\mathscr{X}$ and let $\mathcal{P}_{\mathbb{X}}$ be the space of probability
measures on $\mathbb{X}$ with the $\sigma$-field
$\mathscr{P}_{\mathbb{X}}$
generated by the topology of weak convergence. If $\alpha$ is a strictly positive finite
measure on $\mathbb{X}$ with total mass $a>0$, $Y$ is a $\mathbb{X}$-valued
random variable (r.v.) distributed according to $\alpha_{0}:=\alpha/a$ and $\theta$ is a r.v.
independent of $Y$ and distributed according to a~Beta distribution with parameter
$(1,a)$ then, Theorem 3.4 in \cite{Set94} 
implies that a Dirichlet process
$P$ on $\mathbb{X}$ with parameter $\alpha$ uniquely satisfies the distributional
equation
%
\begin{equation}\label{ftchain1}
P\stackrel{\mathrm{d}}{=}\theta\delta_{Y}+(1-\theta)P,
\end{equation}
where all the random elements on the right-hand side of \eqref{ftchain1} are
independent. All the r.v.s introduced in this paper are meant to be assigned on a
probability space $(\Omega,\mathscr{F},\mathbb{P})$ unless otherwise stated.
In \cite{Fei89}, 
\eqref{ftchain1} is recognized as the distributional equation
for the unique invariant measure of a measure-valued Markov chain $\{P_{m},m\geq0\}$
defined via the recursive identity
%
\begin{equation}\label{ftchain2}
P_{m}=\theta_{m}\delta_{Y_{m}}+(1-\theta_{m})P_{m-1},\qquad m\geq1,
\end{equation}
where $P_{0}\in\mathcal{P}_{\mathbb{X}}$ is arbitrary, $\{Y_{m},m\geq1\}$ is a sequence
of $\mathbb{X}$-valued r.v.s independent and identically distributed as $Y$ and
$\{\theta_{m},m\geq1\}$ is a sequence of r.v.s, independent and identically distributed as
$\theta$ and independent of $\{Y_{m},m\geq1\}$. We term $\{P_{m},m\geq0\}$ as the Feigin--Tweedie
Markov chain. By investigating the functional Markov chain \mbox{$\{G_{m},m\geq0\}$}, with
$G_{m}:=\int_{\mathbb{X}}g(x)P_{m}(\mathrm{d}x)$ for any $m\ge0$ and for any measurable linear
function $g\dvtx\mathbb{X}\mapsto\mathbb{R}$, \cite{Fei89} 
provide
properties of the corresponding linear functional of a Dirichlet process. In particular,
the existence of the linear functional $G:=\int_{\mathbb{X}}g(x)P(\mathrm{d}x)$ of the Dirichlet
process $P$ is characterized according to the condition $\int_{\mathbb{X}}\log(1+|g(x)|)\alpha(\mathrm{d}x)<+\infty$;
these functionals were considered by \cite{Han81} 
and their existence was also investigated by \cite{Dos82} 
who referred
to them as moments, as well as by \cite{Yam84} 
and \cite{Cif90}. 
Further developments of the linear functional Markov chain $\{G_{m},m\geq0\}$ are provided
by \cite{Gug01,Jar02} 
and more recently by \cite{Erh08}. 

Starting from the distributional equation \eqref{ftchain1}, a constructive
definition of the Dirichlet process has been proposed by \cite{Set94}. 
If $P$
is a Dirichlet process on $\mathbb{X}$ with parameter $\alpha=a\alpha_0$, then
$P=\sum_{1\leq i\leq\infty }p_{i}\delta_{Y_{i}}$ where $\{Y_{i},i\geq1\}$ is a sequence of independent r.v.s
identically distributed according to $\alpha_{0}$ and $\{p_{i},i\geq1\}$ is a sequence of r.v.s
independent of $\{Y_{i},i\geq1\}$ and derived by the so-called stick breaking construction,
that is, $p_{1}=w_{1}$ and $p_{i}=w_{i}\prod_{1\leq j\leq i-1}(1-w_{j})$ for $i>1$, with $\{w_{i},i\geq1\}$
being a sequence of independent r.v.s identically distributed according to a Beta distribution
with parameter $(1,a)$. Then, equation \eqref{ftchain1} arises by considering
\[
P=p_{1}\delta_{Y_{1}}+(1-w_{1})\sum_{i=2}^{\infty}\tilde{p}_{i}\delta_{Y_{i}},
\]
where now $\tilde{p}_{2}=w_{2}$ and $\tilde{p}_{i}=w_{i}\prod_{2\leq j\leq i-1}(1-w_{j})$ for $i>2$.
Thus, it is easy to see that $\tilde{P}:=\sum_{2\le i\leq\infty}\tilde{p}_{i}\delta_{Y_{i}}$ is also a Dirichlet
process on $\mathbb{X}$ with parameter $\alpha$ and it is independent of the pairs of r.v.s
$(p_{1},Y_{1})$. If we would extend this idea to $n$ initial samples, we should consider
writing
\[
P=\theta\sum_{i=1}^{n}\biggl(\frac{p_{i}}{\theta}\biggr)\delta_{Y_{i}}+(1-\theta)\tilde{P},
\]
where $\theta=\sum_{1\leq i\leq n}p_{i}=1-\prod_{1\leq i\leq n}(1-w_{i})$ and $\tilde{P}$ is a Dirichlet
process on $\mathbb{X}$ with parameter $\alpha$ independent of the random vectors
$(p_{1},\ldots,p_{n})$ and $(Y_{1},\ldots,Y_{n})$. However, this is not an easy
extension since the distribution of $\theta$ is unclear, and moreover $\theta$ and
$\sum_{1\leq i\leq n}(p_{i}/\theta)\delta_{Y_{i}}$ are not independent. For this reason, in \cite{Fav07}
an alternative distributional equation has been
introduced. Let $\alpha$ be a strictly positive finite measure on $\mathbb{X}$ with
total mass $a>0$ and let $\{Y_{j},j\geq1\}$ be a $\mathbb{X}$-valued P\'olya sequence with parameter
$\alpha$ (see \cite{Bla73}), that is, $\{Y_{j},j\geq1\}$ is a sequence of $\mathbb{X}$-valued r.v.s
characterized by the following
predictive distributions
\[
\P(Y_{j+1}\in A|Y_{1},\ldots,Y_{j})=\frac{1}{a+j}\alpha(A)+\frac{1}{a+j}\sum_{i=1}^{j}\delta_{Y_{i}}(A),\qquad j\geq 1,
\]
and $\P(Y_1\in A)=\alpha(A)/a$, for any $A\in\mathscr{X}$. The sequence $\{Y_{j},j\geq1\}$ is exchangeable, that is,
for any $j\geq1$ and any permutation $\sigma$ of the indexes $(1,\ldots,j)$, the
law of the r.v.s $(Y_{1},\ldots,Y_{j})$ and $(Y_{\sigma(1)},\ldots,Y_{\sigma(j)})$ coincide; in particular, according
to the celebrated de Finetti representation theorem, the P\'olya sequence $\{Y_{j},j\geq1\}$
is characterized by a so-called de Finetti measure, which is the law of a Dirichlet
process on $\mathbb{X}$ with parameter $\alpha$. For a fixed integer
$n\geq1$, let $(q_{1}^{(n)},\ldots,q_{n}^{(n)})$ be a random vector distributed according
to the Dirichlet distribution with parameter $(1,\ldots,1)$, $\sum_{1\leq i\leq n }q_i^{(n)}=1$,
and let $\theta$ be a r.v. distributed according to a Beta distribution with parameter
$(n,a)$ such that $\{Y_{i},i\geq1\}$, $(q^{(n)}_{1},\ldots,q^{(n)}_{n})$ and $\theta$ are mutually
independent. Moving from such a collection of random elements, Lemma 1 in
\cite{Fav07} 
implies that a Dirichlet process $P^{(n)}$ on $\mathbb{X}$
with parameter $\alpha$ uniquely satisfy the distributional equation
%
\begin{equation}\label{ftchain3}
P^{(n)}\stackrel{\mathrm{d}}{=}\theta\sum_{i=1}^{n}q_{i}^{(n)}
\delta_{Y_{i}}+(1-\theta)P^{(n)},
\end{equation}
where all the random elements on the right-hand side of \eqref{ftchain3} are
independent. In order to emphasize the additional parameter $n$, we used an
upper-script $(n)$ on the Dirichlet process $P$ and on the random vector
$(q_1^{(n)},\ldots,q_{n}^{(n)})$. It can be easily checked that equation~\eqref{ftchain3}
generalizes~\eqref{ftchain1}, which can be recovered by setting $n=1$.

In the present paper, our aim is to further investigate the distributional
equation \eqref{ftchain3} and its implications in Bayesian nonparametrics theory
and methods. The first part of the paper is devoted to investigate the random element
$\sum_{1\leq i\leq n }q_{i}^{(n)}\delta_{Y_{i}}$ in \eqref{ftchain3} which is recognized to be
the random probability measure (r.p.m.) at the $n$th step of a measure-valued
Markov chain defined via the recursive identity
%
\begin{equation}\label{nuovacatena}
\sum_{i=1}^{n}q_{i}^{(n)}\delta_{Y_{i}}=W_{n}\delta_{Y_{n}}+(1-W_{n})\sum_{i=1}^{n-1}q_{i}^{(n-1)}\delta_{Y_{i}},\qquad
n\geq1,
\end{equation}
where $\{W_{n},n\geq1\}$ is a~sequence of independent r.v.s, each $W_n$ distributed
according a~Beta distribution with parameter $(1,n-1)$, $q_{i}^{(n)}=W_{i}\prod_{i+1\leq j\le n}(1-W_{j})$
for $i=1,\ldots,n$ and $n\geq1$ and the sequence $\{W_{n},n\geq1\}$ is independent from
$\{Y_{n},n\geq1\}$. More generally, we observe that the measure-valued Markov chain
defined via the recursive identity \eqref{nuovacatena} can be extended by considering,
instead of a P\'olya sequence $\{Y_{n},n\geq1\}$ with parameter~$\alpha$,
any exchangeable sequence $\{Z_{n},n\geq1\}$ characterized by some de Finetti measure on~$\mathcal{P}_{\mathbb{X}}$ and such that $\{W_{n},n\geq1\}$ is independent from
$\{Z_{n},n\geq1\}$.
Asymptotic properties for this new class of measure-valued Markov chains are derived and
some linkages to Bayesian nonparametric mixture modelling are discussed.
In particular, we remark how it is closely related to a well-known recursive algorithm
introduced in  \cite{Newt99} for estimating the underlying mixing distribution in
mixture models, the so-called Newton's algorithm.

In the second part of the paper, by using finite and asymptotic properties of the
r.p.m. $\sum_{1\leq i\leq n }q_{i}^{(n)}\delta_{Y_{i}}$ and by following the original idea of \cite{Fei89}, 
we define and investigate from \eqref{ftchain3}
a class of measure-valued Markov chain $\{P^{(n)}_{m},m\geq0\}$ which generalizes the
Feigin--Tweedie Markov chain, introducing a fixed integer parameter $n$.
Our aim is in providing features of the Markov chain
$\{P^{(n)}_{m},m\geq0\}$ in order to verify if it preserves some of the
properties characterizing the Feigin--Tweedie Markov chain; furthermore,
we are interested in analyzing asymptotic (as $m$ goes to $+\infty$) properties of the associated linear functional Markov chain
$\{G^{(n)}_{m},m\geq0\}$ with $G^{(n)}_{m}:=\int_{\mathbb{X}}g(x)P^{(n)}_{m}(\mathrm{d}x)$ for any
$m\ge0$ and for any function $g\dvtx\mathbb{X}\mapsto\mathbb{R}$ such that $\int_{\mathbb{R}}\log(1+|g(x)|)\alpha(\mathrm{d} x)<+\infty$.
In particular, we show that the Feigin--Tweedie Markov chain
$\{P_{m},m\geq0\}$ sits in a larger class of measure-valued Markov chains $\{P^{(n)}_{m},m\geq0\}$
parametrized by an integer number $n$ and still having the law of a Dirichlet process
with parameter $\alpha$ as unique invariant measure. The role of the further
parameter $n$ is discussed in terms of new potential applications of the Markov chain
$\{P^{(n)}_{m},m\geq0\}$ with respect to the the known applications of the Feigin--Tweedie
Markov chain.

Following these guidelines, in Section \ref{sec2} we introduce a new class of
measure-valued Markov chains $\{Q_{n},n\geq1\}$ defined via exchangeable sequences of
r.v.s; asymptotic results for $\{Q_{n},n\geq1\}$ are derived and applications related to
Bayesian nonparametric mixture modelling
are discussed. In Section \ref{secthree}, we show that the Feigin--Tweedie Markov chain $\{P_{m},m\geq0\}$ sits in a larger
class of measure-valued Markov chains $\{P^{(n)}_{m},m\geq0\}$, which is investigated in comparison with $\{P_{m},m\geq0\}$.
In Section \ref{sec4}, some concluding remarks and future
research lines are presented.

\section{A class of measure-valued Markov chains and Newton's algorithm}\label{sec2}
Let $\{W_{n},n\geq1\}$ be a sequence of independent r.v.s such that $W_{1}=1$ almost surely and~$W_n$
has Beta distribution with parameter $(1,n-1)$ for
$n\geq 2$. Moreover, let \mbox{$\{Z_{n},n\geq1\}$} be a sequence of $\mathbb{X}$-valued
exchangeable r.v.s independent from $\{W_{n},n\geq1\}$ and characterized by some de Finetti measure on
$\mathcal{P}_{\mathbb{X}}$.
Let us consider the measure-valued Markov chain $\{Q_{n},n\geq1\}$ defined via the recursive identity
%
\begin{equation}\label{eqPrec}
Q_{n}=W_{n}\delta_{Z_{n}}+(1-W_{n})Q_{n-1},\qquad n\geq 1.
\end{equation}
In the next theorem, we provide an alternative representation of $Q_{n}$ and show that
$Q_n(\omega)$ converges weakly to some limit probability $Q(\omega)$ for almost
all $\omega\in\Omega$, \textit{that is}, for each $\omega$ in some set $A\in\mathscr{F}$
with $\mathbb{P}(A)=1$. In short, we use notation $Q_{n}\Rightarrow Q$ a.s.-$\mathbb{P}$.

\begin{theorem}\label{asconv}
Let $\{Q_{n},n\geq 1\}$ be the Markov chain defined by \eqref{eqPrec}. Then:

\begin{enumerate}[(ii)]
\item[(i)] an equivalent representation of $Q_n$, $n=1,2,\ldots$ is
%
\begin{equation}\label{eqqrepr}
Q_{n}=\sum_{i=1}^{n}q_i^{(n)}\delta_{Z_{i}},
\end{equation}
where $\sum_{1\le i\le n }q_i^{(n)}=1$, $q^{(n)}=(q^{(n)}_1,\ldots,q^{(n)}_n )$
has Dirichlet distribution with parameter $(1,1,\ldots,1)$, and $\{q^{(n)},n\geq 1 \}$ and
$\{Z_{n},n\geq1\}$ are independent.

\item[(ii)] There exists a r.p.m. $Q$ on $(\Omega,\mathscr{F},\mathbb{P})$ such that, as
$n\rightarrow+\infty$,
\[
Q_{n}\Rightarrow Q,\qquad\mbox{\textup{a.s.}-}\mathbb{P},
\]
where the law of $Q$ is the de Finetti measure of the sequence $\{Z_{n},n\geq1\}$.
\end{enumerate}
\end{theorem}

\begin{pf}
As far as (i) is concerned, by repeated application of the recursive identity \eqref{eqPrec},
it can be checked that, for any $n\geq1$,
\[
Q_{n}=\sum_{i=1}^{n}W_{i}\prod_{j=i+1}^{n}(1-W_j)\delta_{Z_{i}},
\]
where $W_1=1$ almost surely and $\prod_{i+1\leq j\leq n}(1-W_j)$ is defined to be 1 when
$i=n$. Defining $q_i^{(n)}:=W_{i}\prod_{i+1\leq j\leq n}(1-W_j)$, $i=1,\ldots,n$,
it is straightforward to show that
$q^{(n)}=(q^{(n)}_1,\ldots,q^{(n)}_n )$ has the
Dirichlet distribution with parameter $(1,1,\ldots,1)$ and $\sum_{1\leq i\leq n }q_i^{(n)}=1$,
so that \eqref{eqqrepr} holds.

Regarding (ii), by the definition of
the Dirichlet distribution, an equivalent representation of
\eqref{eqqrepr} is\vspace*{2pt}
\[
Q_{n}=\sum_{i=1}^{n}q_{i}^{(n)}\delta_{Z_{i}}=
\sum_{i=1}^{n}\frac{\lambda_{i}}{\sum_{j=1}^{n}\lambda_{j}}\delta_{Z_{i}},
\]
where $\{\lambda_n, n\geq 1\}$ is a sequence of r.v.s independent and identically distributed according to
standard exponential distribution, independent from  $\{Z_{n},n\geq1\}$.
Let $g\dvtx\mathbb{X}\rightarrow\mathbb{R}$ be any bounded continuous function, and consider\vspace*{2pt}
\[
G_n=\int_{\mathbb X}g\,\mathrm{d}Q_n=\sum_{i=1}^n\frac{\lambda_i}{\sum_1^n\lambda_j} g(Z_i)=
\frac{{\sum_{i=1}^n\lambda_ig(Z_i)}/n}{{\sum_1^n\lambda_i}/n}.
\]
The expression in the denominator converges almost surely to 1 by the strong law of
large numbers.
As far as the numerator is concerned, let $Q$ be the r.p.m. defined on
$(\Omega,\mathscr{F},\mathbb{P})$, such that the r.v.s $\{Z_n, n\geq 1\}$ are
independent and identically distributed conditionally on $Q$;
the existence of such a random element is guaranteed by the de Finetti representation theorem
(see, e.g., \cite{Sche95}, 
Theorem 1.49). It can be shown that
$\{\lambda_n g(Z_n), n\geq 1\}$ is a sequence of exchangeable r.v.s and, if $t_1,\ldots,t_n\in\R$,\vspace*{2pt}
\begin{eqnarray*}
&&\mathbb{P}\bigl(\lambda_1 g(Z_1)\leq t_1,\ldots,\lambda_n g(Z_n) \leq t_n\bigr) \\
&&\quad=\int_{(0,+\infty)^n}\mathbb{P}\bigl(\lambda_1g(Z_1)\leq t_1,\ldots,\lambda_ng(Z_n) \leq t_n|
\lambda_1,\ldots,\lambda_n\bigr)\prod_{i=1}^n \e^{-\lambda_i}\,\mathrm{d}\lambda_i \\
&&\quad =\int_{(0,+\infty)^n}\int_{{\mathcal{P}}_{\R}}\prod_{i=1}^n F_{Q^*}\biggl(\frac{t_{i}}{\lambda_{i}}\biggr)
\mu(\mathrm{d}Q^*)\prod_{i=1}^n \e^{-\lambda_i}\,\mathrm{d}\lambda_i\\
&&\quad=\int_{\mathcal{P}_{\R}}\prod_{i=1}^n\biggl(\int_0^{+\infty}F_{Q^*}\biggl(\frac{t_i}{\lambda_i}\biggr)
\e^{-\lambda_i}\,\mathrm{d}\lambda_i\biggr)\mu(\mathrm{d}Q^*),
\end{eqnarray*}
where $Q^*(A,\omega):=P(g^{-1}(A),\omega)$, $\omega\in\Omega$, $A\in\mathscr{R}$, is a r.p.m.
with trajectories in ${\mathcal{P}_{\R}}$, and~$F_{Q^*}$ denotes the random distribution
relative to $Q^*$. This means that, conditionally on $Q^*$,
$\{\lambda_n g(Z_n), n\geq 1\}$ is a sequence of r.v.s independent and identically
distributed according to the random distribution (evaluated in $t$)
\[
\int_0^{+\infty}F_{Q^*}\biggl(\frac{t}{y}\biggr)\e^{-y}\,\mathrm{d}y.
\]
Of course, $\E|\lambda_1 g(Z_1)|=\E(\lambda_1)\E|g(Z_1)|<+\infty$ since $g$ is bounded.
As in \cite{ChowTei97},
Example~7.3.1, this condition implies
\begin{eqnarray*}
\frac{1}{n} \sum_{i=1}^n \lambda_i g(Z_i)&\stackrel{\mathrm{a.s.}}\rightarrow&
\E(\lambda_1
g(Z_1)|Q^*)=\int_{\R}t\,\mathrm{d}\biggl(\int_0^{+\infty}F_{Q^*}\biggl(\frac{t}{y}\biggr)\e^{-y}\,\mathrm{d}y\biggr)\\
&=& \int_{\R} u Q^*(\mathrm{d}u)=\int_{\mathbb{X}}g(x) Q(\mathrm{d}x),
\end{eqnarray*}
so that $G_n\rightarrow \int_{\mathbb{X}}g(x) Q(\mathrm{d}x)$ a.s.-$\mathbb{P}$.
By Theorem 2.2 in \cite{Ber06}, 
it follows that $Q_{n}\Rightarrow Q$
a.s.-$\mathbb{P}$  as $n\rt +\infty$.
\end{pf}

Throughout the paper, $\alpha$ denotes a strictly positive and finite measure on $\mathbb{X}$
with total
mass~$a$, unless otherwise stated. If the exchangeable sequence $\{Z_{n},n\geq1\}$
is the 
P\'olya sequence with parameter $\alpha$,
then by Theorem \ref{asconv}(i) $\{Q_{n},n\geq1\}$ is the Markov chain defined via the
recursive identity \eqref{nuovacatena}; in particular, by Theorem \ref{asconv}(ii),
$Q_{n}\Rightarrow Q$ $\mathrm{a.s.\mbox{-}}\mathbb{P}$ where $Q$ is a Dirichlet process on
$\mathbb{X}$ with parameter $\alpha$.
This means that, for any fixed integer $n\geq1$, the r.p.m.
$Q_{n}$ can be interpreted as an
approximation of a Dirichlet process with parameter $\alpha$.
In Appendix \ref{appa1}, we present an alternative proof of the weak convergence
(convergence of the finite dimensional distribution) of $\{Q_{n},n\geq1\}$ to a
Dirichlet process on $\mathbb{X}$ with
parameter $\alpha$, using a combinatorial technique.
As a byproduct of this proof, we obtain an explicit expression for the
moment of order $(r_{1},\ldots,r_{k}$) of the $k$-dimensional P\'olya distribution.

A straightforward generalization of the Markov chain $\{Q_{n},n\geq0\}$ can be obtained
by considering a nonparametric hierarchical mixture model.
Let $k\dvtx\mathbb{X}\times\Theta\rightarrow\mathbb{R}^{+}$ be a kernel, that is,
$k(x,\vartheta)$ is a measurable function such that
$x\mapsto k(x,\vartheta)$ is a density with respect to some $\sigma$-finite measure
$\lambda$ on $\mathbb{X}$, for any fixed $\vartheta\in\Theta$, where $\Theta$
is a
Polish space (with the usual Borel $\sigma$-field). Let
$\{Q_{n},n\geq1\}$ be the Markov chain defined via \eqref{eqPrec}. Then
for each $x\in\mathbb{X}$ we introduce a real-valued Markov chain
$\{f^{(Q)}_{n}(x),n\ge1\}$ defined via the recursive identity
%
\begin{equation}\label{catenamistura}
f^{(Q)}_{n}(x)=W_{n}k(x,\vartheta_{n})+(1-W_{n})f^{(Q)}_{n-1}(x),\qquad n\geq 1,
\end{equation}
where
\[\label{eqfn}
f^{(Q)}_{n}(x)=\int_{\Theta}k(x,\vartheta)Q_{n}(\mathrm{d}\vartheta).
\]
By a straightforward application of Theorem 2.2 in  \cite{Ber06},
for any fixed $x\in\mathbb{X}$, when
$\vartheta\mapsto k(x,\vartheta)$ is continuous for all  $x\in\mathbb{X}$ and bounded by a function $h(x)$,
as $n\rightarrow+\infty$, then\vspace*{2pt}
%
\begin{equation}\label{convdens}
f^{(Q)}_{n}(x)\rightarrow f^{(Q)}(x):=\int_{\Theta}k(x,\vartheta)Q(\mathrm{d}\vartheta),\qquad\mathrm{a.s.\mbox{-}}\mathbb{P},
\end{equation}
where $Q$ is the limit r.p.m. in Theorem \ref{asconv}. For instance, if $Q$ is a
Dirichlet process on $\mathbb{X}$ with parameter $\alpha$, $f^{(Q)}$ is precisely the
density in the Dirichlet process mixture model introduced by \cite{Lo84}. 
When $h(x)$ is a $\lambda$-integrable function, not
only the limit $f^Q(x)$ is a~random density,
but a stronger result than \eqref{convdens} is achieved.

\begin{theorem}\label{rigo}
If $\vartheta\mapsto k(x,\vartheta)$ is continuous for all  $x\in\mathbb{X}$ and bounded by a $\lambda$-integrable function
$h(x)$, then\vspace*{2pt}
\[
\lim_{n\rightarrow+\infty}\int_{\mathbb{X}}
\bigl|f^{(Q)}_{n}(x)-f^{(Q)}(x)\bigr|\lambda(\mathrm{d}x)\rightarrow 0,\qquad\mathrm{a.s.\mbox{-}}\mathbb{P},
\]
where $Q$ is the limit r.p.m. in Theorem \ref{asconv}.
\end{theorem}

\begin{pf}
The functions $f^{(Q)}_n$ and  $f^{(Q)}(x)=\int_{\Theta}k(x,\vartheta)Q(\mathrm{d}\vartheta)$,
defined
on $\mathbb{X}\times \Omega$, are  $\mathscr{X}\otimes\mathscr{F}$-measurable, by a monotone class argument.
In fact, by kernel's definition, $(x,\vartheta)\mapsto k(x,\vartheta)$ is
$\mathscr{X}\otimes\calB(\Theta)$-measurable. Moreover, if $k=\mathbbl{1}_A\mathbbl{1}_B$, $A\in \mathscr{X}$
and $B\in\calB(\Theta)$, then\vspace*{2pt}
\[
f^{(Q)}(x,\omega)=\int k(x,\vartheta)Q(\mathrm{d}\vartheta;\omega)=\mathbbl{1}_A(x)Q(B;\omega)
\]
is $\mathscr{X}\otimes\mathscr{F}$-measurable. Let
${\mathcal C}=\{C\in \mathscr{X}\otimes\calB(\Theta)\dvt
\int \mathbbl{1}_C(x,\vartheta) Q(\mathrm{d}\vartheta;\omega)\mbox{ is } \mathscr{X}\otimes\mathscr{F}\mathrm{\mbox{-}measurable}\}$.
Since $\mathcal C$ contains the rectangles, it contains the field generated by
rectangles, and, since~$\mathcal C$ is a monotone class,
$\mathcal C=\mathscr{X}\otimes\calB(\Theta)$. The assertion holds for $f^{(Q)}$ of the form\vspace*{2pt}
\[
f^{(Q)}(x)=\int_{\Theta}k(x,\vartheta)Q(\mathrm{d}\vartheta)
\]
since there exist a sequence of simple function on rectangles which
converges pointwise to~$k$.
Therefore,
$A:=\{(\omega,x))\dvt f^{(Q)}_n(\omega,x)$ does not converge to $f^{(Q)}(\omega,x)\}\in\mathscr{F}\otimes\mathscr{X}$.
Then, by Fubini's theorem,\vspace*{2pt}
\begin{eqnarray*}
&&\lefteqn{\int\lambda\bigl\{x\dvt f^{(Q)}_n(\omega,x)\textrm{ does not converge to } f^{(Q)}(\omega,x)\bigr\} \mathbb{P}(\mathrm{d}\omega)}\\
&&\quad= \int\int
\mathbbl{1}_A(\omega,x)\lambda(\mathrm{d}x)\mathbb{P}(\mathrm{d}\omega)\\
&&\quad=\int\mathbb{P}\bigl\{\omega\dvt f^{(Q)}_n(\omega,x)\textrm{ does not converge to } f^{(Q)}(\omega,x) \bigr\} \lambda(\mathrm{d}x)\\
&&\quad=\int 0\lambda(\mathrm{d}x)=0.
\end{eqnarray*}
Hence, $\mathbb{P}(H)=1$ where $H$ is the set of $\omega$ such that
$\lambda\{x\dvt f^{(Q)}_n(\omega,x)$
does not converge to  $f^{(Q)}(\omega,x)\}=0$.
For any $\omega$ fixed in $H$, it holds $f^{(Q)}_n(\omega,\cdot)\rightarrow f^{(Q)}(\omega,\cdot)$, $\lambda$-a.e., so that by the
Scheff\'e's theorem we have
\[
\lim_{n\rightarrow+\infty}\int_{\mathbb{X}}\bigl|f^{(Q)}_{n}(\omega,x)-f^{(Q)}(\omega,x)\bigr|\lambda(\mathrm{d}x)\rightarrow 0.
\]
The theorem follows since $\mathbb{P}(H)=1$.
\end{pf}

We conclude this section by remarking an interesting linkage between the Markov chain
$\{Q_{n},n\geq1\}$ and the so-called Newton's algorithm, originally introduced in
\cite{Newt99} for estimating the mixing density when a finite sample is available
from the corresponding mixture model.
See also See also \cite{Newt98} 
and \cite{Newt02}. 
Briefly, suppose that $X_{1},\ldots,X_{n}$ are $n$ r.v.s independent and identically
distributed according to the density function
%
\begin{equation}\label{mixmod}
\widetilde{f}(x)=\int_{\Theta}k(x,\vartheta)\widetilde{Q}(\mathrm{d}\vartheta),
\end{equation}
where $k(x,\vartheta)$ is a known kernel dominated by a $\sigma$-finite measure $\lambda$
on $\mathbb{X}$; assume that the mixing distribution $\widetilde{Q}$ is absolutely continuous
with respect to some
$\sigma$-finite measure $\mu$ on $\Theta$. \cite{Newt02} 
proposed to
estimate $\widetilde{q}=\mathrm{d}\widetilde{Q}/\mathrm{d}\mu$ as follows: fix an initial estimate
$\hat{q}_{1}$ and a~sequence of weights $w_{1},w_{2},\ldots,w_{n}\in(0,1)$. Given $X_{1},\ldots,X_{n}$
independent and identically distributed observations
from $\widetilde{f}$, compute
\[\label{alg}
\hat q_{i}(\vartheta)=(1-w_{i})\hat q_{i-1}(\vartheta)+w_{i}\frac{k(x_{i},\vartheta)
\hat q_{i-1}(\vartheta)}{\int_{\Theta}k(x_{i},\vartheta)\hat q_{i-1}(\vartheta)\mu(\mathrm{d}\vartheta)},\qquad\vartheta\in\Theta
\]
for $i=2,3,\ldots,n$ and produce $\hat q_{n}$ as the final estimate.
We refer to \cite{Gho06,Mar08,Tokdar08},  
and~\cite{MarTok09} 
for a recent
wider investigation of the Newton's algorithm.
Here we show how the Newton's algorithm is connected to the
measure-valued Markov chain $\{Q_{n},n\geq1\}$.

Let us consider $n$ observations from the nonparametric hierarchical mixture model,
that is, $X_i|\vartheta_i\sim k(\cdot,\vartheta_i)$ and
$\vartheta_i|Q\sim Q$ where $Q$ is a r.p.m. If we observed
$\{\vartheta_{i},i\geq1\}$, then by virtue of (ii) in Theorem \ref{asconv}, we could
construct a sequence of distributions
\[
Q_i=W_i\delta_{\vartheta_i}+(1-W_i)Q_{i-1},\qquad i=1,\ldots,n
\]
for estimating the limit r.p.m. $Q$, where $\{W_{i},i\geq1\}$ is a sequence of independent r.v.s such that
$W_{1}=1$ almost surely and $W_{i}$ has
Beta distribution with parameters $(1,i-1)$.
This approximating sequence is precisely the sequence (\ref{eqPrec}).
Therefore, taking the expectation
of both sides of the previous recursive equation, 
and defining $\widetilde{Q}_{i}:=\E[Q_{i}]$, $w_{i}=\E[W_{i}]=1/i$, we have
%
\begin{equation}\label{eqQtilde}
\widetilde{Q}_{i}=w_{i}\delta_{\vartheta_{i}}+(1-w_{i})\widetilde{Q}_{i-1},\qquad
i=1,\ldots,n,
\end{equation}
which can represent a predictive distribution for $\vartheta_{i+1}$, and hence an estimate for $Q$.

However, instead of observing the sequence $\{\vartheta_i,i\geq1\}$, it is actually the sequence $\{X_i,i\geq1\}$ which
is observed; in particular, we can assume that $X_{1},\ldots,X_{n}$
are $n$ r.v.s independent and identically distributed according to the density function \eqref{mixmod}.
Therefore, instead of \eqref{eqQtilde}, we consider
\[
\widetilde{Q}_i(\mathrm{d}\vartheta)=(1-w_i)\widetilde{Q}_{i-1}(\mathrm{d}\vartheta)+ w_i
\frac{k(x_i,\vartheta)\widetilde{Q}_{i-1}(\mathrm{d}\vartheta)}{\int_{\Theta} k(x_i,\vartheta)\widetilde{Q}_{i-1}(\mathrm{d}\vartheta)},\qquad
i=1,\ldots,n,
\]
where $\delta_{\vartheta_i}$ in \eqref{eqQtilde} has been substituted (or estimated,
if you prefer) by
$ {k(x_i,\vartheta)\widetilde{Q}_{i-1}(\mathrm{d}\vartheta)}/\break{\int_{\Theta} k(x_i,\vartheta)\widetilde{Q}_{i-1}(\mathrm{d}\vartheta)}$.
Finally, observe that, if
$\widetilde Q_{i}$ is absolutely continuous, with respect to some $\sigma$-finite
measure~$\mu$ on $\Theta$, with density $\tilde q_i$,
for $i=1,\ldots,n$, then we can write
%
\begin{equation}\label{Newt}
\widetilde{q}_i(\vartheta)=(1-w_{i})\widetilde{q}_{i-1}(\vartheta)+w_{i}
\frac{k(x_i,\vartheta)\widetilde{q}_{i-1}(\vartheta)}
{\int_{\Theta} k(x_i,\vartheta)\widetilde{q}_{i-1}(\vartheta)\mu(\mathrm{d}\vartheta)},\qquad i=1,\ldots,n,
\end{equation}
which is precisely a recursive estimator of a mixing distribution proposed by
\cite{Newt02} 
when the weights are fixed to be
$w_{i}=1/i$ for $i=1,\ldots,n$ and the initial estimate is $\E[\delta_{\vartheta_{1}}]$.

\section{A generalized Feigin--Tweedie Markov chain} \label{secthree}

In this section our aim is to define and investigate a class of measure-valued Markov chain which
generalizes the Feigin--Tweedie Markov chain introducing a fixed integer parameter~$n$, and still has the law of a Dirichlet process
with
parameter $\alpha$ as the unique invariant measure.
The starting point is the
distributional equation
(\ref{ftchain3}) introduced by~\cite{Fav07}; 
see Appendix \ref{appa2} for an alternative proof
of the solution of the distributional equation~\eqref{ftchain3}.
All the proofs of Theorems in this section are in Appendix \ref{appa3} for the ease of reading.

For a fixed integer $n\geq1$, let $\theta:=\{\theta_{m},m\geq1\}$ be a sequence of
independent r.v.s with Beta
distribution with parameter $(n,a)$, $q^{(n)}:=\{(q_{m,1}^{(n)},\ldots,q_{m,n}^{(n)}),m\geq1\}$, with
$\sum_{1\leq i\leq n} q_{m,i}^{(n)}=1$ for any $m>0$, be a sequence of independent r.v.s
identically
distributed according to a Dirichlet distribution with parameter $(1,\ldots,1)$ and
$Y:=\{(Y_{m,1},\ldots,Y_{m,n})$, $m\geq 1\}$ be sequence of independent r.v.s from a
P\'olya sequence with parameter $\alpha$. Moving from such collection of random elements,
 for each
fixed integer $n\geq1$ we define the measure-valued Markov chain $\{P^{(n)}_{m},m\geq0\}$ via the recursive identity
%
\begin{equation}\label{ftchain4}
P^{(n)}_{m}=\theta_{m}\sum_{i=1}^{n}q_{m,i}^{(n)}\delta_{Y_{m,i}}+(1-\theta_{m})P^{(n)}_{m-1},\qquad
m\geq1,
\end{equation}
where $P^{(n)}_{0}\in\mathcal{P}_{\mathbb{X}}$ is arbitrary.
By construction, the Markov chain $\{P_{m},m\geq0\}$
proposed by \cite{Fei89} 
and defined via the recursive identity
\eqref{ftchain2} can be recovered from $\{P^{(n)}_{m},m\geq0\}$ by setting $n=1$.
Following the original idea of \cite{Fei89}, 
by equation \eqref{ftchain4} we have defined the Markov chain
$\{P^{(n)}_{m},m\geq0\}$ from a distributional equation having
as the unique solution the Dirichlet process. In particular, the Markov chain $\{P^{(n)}_{m},m\geq0\}$ is defined
from the distributional equation \eqref{ftchain3} which generalizes \eqref{ftchain1} substituting  the random
probability measure $\delta_{Y}$ with the random
convex linear combination $\sum_{1\leq i\leq n}q_{i}^{(n)}\delta_{Y_{i}}$, for any fixed positive integer $n$.
Observe that$\sum_{1\leq i\leq n}q_{i}^{(n)}\delta_{Y_{i}}$ is an example of the r.p.m. $Q_n$
defined in \eqref{eqqrepr}  and investigated in the previous section, when $\{Z_i \}$ is
given by the P\'olya sequence $\{Y_i\}$ with parameter $\alpha$. In particular,
Theorem \ref{asconv} shows that $Q_n$ a.s.-converges to the Dirichlet process $P$ when
$n$ goes to infinity; however here we assume a different perspective, that is, $n$ is fixed.

As for the case $n=1$, the following result holds.

\begin{theorem}\label{thminvariantPn}
The Markov chain $\{P^{(n)}_{m},m\geq0\}$ has a unique
invariant measure $\Pi$ which is the law of a Dirichlet process $P$ with parameter
$\alpha$.
\end{theorem}

Another property which still holds in the more general case when $n\geq 1$ is the Harris
ergodicity of the functional Markov chain $\{G^{(n)}_{m},m\geq0\}$,
under assumption \eqref{eqcondlog} below. This condition is equivalent to the finiteness of the r.v.
$\int_{\mathbb{X}}|g(x)|P(\mathrm{d}x)$; see also \cite{Cif90}. 

\begin{theorem}\label{thmharrisfunct}
Let $g\dvtx\mathbb{X}\mapsto\mathbb{R}$ be any measurable
function. If
%
\begin{equation}\label{eqcondlog}
\int_{\mathbb{X}}\log\bigl(1+|g(x)|\bigr)\alpha(\mathrm{d}x)<+\infty,
\end{equation}
then the Markov chain $\{G^{(n)}_{m},m\geq0\}$ is Harris ergodic with unique invariant
measure~$\Pi_{g}$, which is the law of the random Dirichlet mean
$\int_{\mathbb{X}}g(x) P(\mathrm{d}x)$.
\end{theorem}

We conclude the analysis of the Markov chain $\{P_{m}^{(n)},m\geq0\}$ by providing some
results on the ergodicity of the
Markov chain $\{G_{m}^{(n)},m\geq0\}$ and by discussing on the rate of convergence.
Let $\mathbb{X}=\mathbb{R}$ and let $\{P^{(n)}_{m},m\geq0\}$ be the Markov chain defined
by \eqref{ftchain4}. In particular, for the rest of the section, we consider the mean
functional Markov chain $\{M^{(n)}_{m},m\geq0\}$ defined recursively by
%
\begin{equation}\label{ftchain11}
M^{(n)}_{m}=\theta_{m}\sum_{i=1}^{n}q_{m,i}^{(n)}Y_{m,i}+(1-\theta_{m})M^{(n)}_{m-1},\qquad
m\geq1,
\end{equation}
where $M^{(n)}_{0}\in\mathbb{R}$ is arbitrary and $n$ is a given positive integer.
From Theorem \ref{thmharrisfunct}, under the condition $\int_{\mathbb{R}}\log(1+|x|)\alpha(\mathrm{d}x)<+\infty$,
the Markov chain $\{M^{(n)}_{m},m\geq0\}$ has the distribution $\mathscr{M}$ of the random
Dirichlet mean $M$ as the unique invariant measure. It is not restrictive to consider
only the chain $\{M^{(n)}_{m},m\geq0\}$, since a more general linear functionals $G$ of a
Dirichlet process on an arbitrary Polish space has the same distribution as the mean
functional of a Dirichlet process with parameter $\alpha_{g}$,
where $\alpha_{g}(B):=\alpha(g^{-1}(B))$ for any $B\in\mathscr{R}$.

\begin{theorem}\label{thmgeoergfunct}
The Markov chain $\{M^{(n)}_{m},m\geq0\}$ satisfies the following properties:
\begin{enumerate}[(iii)]
\item[(i)] $\{M^{(n)}_{m},m\geq0\}$ is a stochastically monotone Markov chain;
\item[(ii)] if further
%
\begin{equation}\label{eqfinitmeanalpha0}
\E[|Y_{1,1}|]=\int_{\mathbb{R}}|x|\alpha_0(\mathrm{d}x)<+\infty,
\end{equation}
then $\{M^{(n)}_{m},m\geq0\}$ is a geometrically ergodic Markov chain;
\item[(iii)] if the support of $\alpha$ is bounded then $\{M^{(n)}_{m},m\geq0\}$ is an uniformly ergodic Markov chain.
\end{enumerate}
\end{theorem}

Recall that the stochastic monotonicity property of $\{M^{(n)}_{m},m\geq0\}$ allows to consider
exact sampling (see \cite{Prop96}) 
for $M$ via
$\{M_m^{(n)},m\geq 0\}$.

\begin{rmk}
\label{remjarner}
Condition \textup{\eqref{eqfinitmeanalpha0}} can be relaxed. If the following condition holds
%
\begin{equation}\label{eqjar}
\E[|Y_{1,1}|^{s}]=\int_{\mathbb{R}}|y|^{s}\alpha_{0}(\mathrm{d}x)<+\infty\qquad\mbox{for some }0<s<1,
\end{equation}
then the Markov chain $\{M_{m}^{(n)},m\geq0\}$ is geometrically ergodic.
See Appendix \textup{\ref{appa3}} for the proof. If, for instance, $\alpha_{0}$ is a Cauchy standard distribution and $a>0$,
condition  \textup{\eqref{eqjar}} is fulfilled so that $\{M_{m}^{(n)},m\geq0\}$ will turn out to
be geometrically ergodic for any fixed integer $n$.
\end{rmk}

From Theorem \ref{asconv},  $T^{(n)}:=\sum_{1\leq i\leq n}q_{i}^{(n)}Y_{i}=\int_\R x
( \sum_{1\leq i\leq n}q_{i}^{(n)}\delta_{Y_{i}})(\mathrm{d}x)$ converges in distribution to the random Dirichlet mean $M$ as $n\rightarrow+\infty$; so it
is clear that, for a fixed integer $n$, the law of $T^{(n)}$ approximates the law of $M$
and that the approximation will be better for $n$ large.
If we reconsider \eqref{ftchain11}, written as
\[
M^{(n)}_{m}=\theta_{m}T_m^{(n)} +(1-\theta_{m})M^{(n)}_{m-1},\qquad m\geq1,
\]
since the innovation term $T^{(n)}$ is an approximation in distribution of the
limit (as $m\rightarrow+\infty$) r.v. $M$, it is intuitive
that the rate of convergence will increase as $n$ gets larger. This is confirmed
by the description of small sets $C^{(n)}$ in
\eqref{ftchain15} (in the proof of Theorem \ref{thmgeoergfunct}).
In fact, under \eqref{eqfinitmeanalpha0} or \eqref{eqjar},
the Markov chain $\{M^{(n)}_{m},m\geq0\}$ is geometrically or uniformly ergodic
since it satisfies a Foster--Lyapunov condition $PV(x):=\int_\R V(y)p(x,\mathrm{d}y)\leq\lambda V(x)+b\mathbbl{1}_{{C}^{(n)}}(x)$ for
a suitable function $V$, a small set $C^{(n)}$ and constants
$b<+\infty$, $0<\lambda<1$. In particular, the small sets $C^{(n)}$ generalize the corresponding small set $C$ obtained in Theorem~1 in \cite{Gug01}
which can be recovered by setting $n=1$, that is, $C=[-K(\lambda),K(\lambda)]$ where
\[
\label{ftchain19}
K(\lambda):=\frac{1-\lambda+1/(1+a)\E[|Y_{1,1}|]}{\lambda-a/(1+a)}.
\]
Here the size of the small set $C^{(n)}$ of $\{M^{(n)}_{m},m\geq0\}$
can be controlled by an additional parameter~$n$, suggesting the upper
bounds of the rate of convergence of the
chain $\{M^{(n)}_{m},m\geq0\}$ depends on $n$ too.

However, if we would establish an explicit upper bound on the rate of convergence,
we would need results like Theorem 2.2 in \cite{Rob00}, 
or Theorems~5.1 and
 5.2 in \cite{Rob99}. 
All these results need a minorization condition to hold for the $m_0$th step transition probability
$p^{(n)}_{m_0}(x,A):=\mathbb{P}(M^{(n)}_{m}\in A|M^{(n)}_{0}=x)$ for any $A\in\mathscr{R}$
and $x\in\mathbb{R}$, for some positive integer~$m_0$ and all $x$ in a small set; in particular, if
$\inf_{x\in C^{(n)}}f(z|x)\geq p_{0}^{(n)}(z)$, where $f(z|x)$ is the density of
$p_{1}^{(n)}(x,\cdot)$ and $p_{0}^{(n)}(z)$ is some
density such that $\varepsilon(n):=\int_{\mathbb{R}}p_{0}^{(n)}(z)\,\mathrm{d}z>0$, then
\[\label{minorcond}
p^{(n)}_{1}(x,A)\geq\varepsilon(n)\int_{A}\frac{p_{0}^{(n)}(z)}{\varepsilon(n)}\,\mathrm{d}z
=\varepsilon(n)\nu(A),\qquad
A\in\mathbb{R}, x\in C^{(n)},
\]
where $\nu$ is a probability measure on $\mathbb{R}$.
In order to check the validity of our intuition that the rate of convergence will
increase as $n$ gets larger,
the function $\varepsilon(n)$ should be increasing with $n$ in order to prove
that the uniform error (when the support of the $Y_{i}$'s is bounded)
in total variation between the law of $M_{m}^{(n)}$ given $M_{0}^{(n)}$
and its limit distribution decreases as $n$ increases. If $f_{T^{(n)}}$ is the density of $T^{(n)}$, which exists since, conditioning on $Y_i$'s,
$T^{(n)}$ is a random Dirichlet mean, then
\[
p^{(n)}_{1}(x,A)=\int_{A}f(z|x)\,\mathrm{d}z
=\int_\R \mathbb{P}\bigl(\theta_{1}y+(1-\theta_{1})x\in A|T^{(n)}=y,M_{0}^{(n)}=x\bigr) f_{T^{(n)}}(y)\,\mathrm{d}y.
\]
Therefore, the density function corresponding to $p^{(n)}_{1}(x,A) $ is
\begin{eqnarray*}
f(z|x) &=& \frac{1}{B(a,n)}\int_\R\frac{(z-x)^{n-1}(y-z)^{a-1}}{(y-x)^{a+n-2}|y-x|}
\mathbbl{1}_{\{(0,1)\}}\biggl(\frac{z-x}{y-x}\biggr)f_{T^{(n)}}(y)\,\mathrm{d}y \\
&=&\cases{\displaystyle
\frac{1}{B(a,n)}\int_{-\infty}^{z}\frac{(x-z)^{n-1}(z-y)^{a-1}}{(x-y)^{a+n-1}}
f_{T^{(n)}}(y)\,\mathrm{d}y, &\quad  if $z < x$, \cr
\displaystyle\frac{1}{B(a,n)}\int_{z}^{+\infty}\frac{(z-x)^{n-1}(y-z)^{a-1}}{(y-x)^{a+n-1}}f_{T^{(n)}}(y)\,\mathrm{d}y,&\quad  if $z>x$}\\
&=&\frac{1}{B(a,n)}\int_{0}^{1}t^{n-2}(1-t)^{a-1}f_{T^{(n)}}\biggl(\frac{z-(1-t)x}{t}\biggr)\,\mathrm{d}t.
\end{eqnarray*}
Unfortunately, the explicit expression of $f_{T^{(n)}}$, which for $n=1$ reduces to the density of $\alpha_{0}$ if it exists,
is not simple; from Proposition~5 in \cite{Reg02} 
for instance,
for $y\in\mathbb{R}$,
\[
f_{T^{(n)}}(y)=\int_{\mathbb{R}^{n}}f_{T^{(n)}}(y;y_{1},\ldots,y_{n})F_{(Y_{1},\ldots,Y_{n})}(\mathrm{d}y_{1},\ldots,y_{n}),
\]
where, when  $y\neq y_{i}$ for $i=1,\ldots,n$,
\begin{eqnarray*}
f_{T^{(n)}}(y;y_{1},\ldots,y_{n})
=\frac{n-1}{\uppi}\int_{0}^{+\infty}\prod_{j=1}^{n}\frac{1}{(1+t^{2}(y_{j}-y)^{2})^{1/2}}
\cos\Biggl(\sum_{j=1}^{n}\arctan\bigl(t(y_{j}-y)\bigr)\Biggr)\,\mathrm{d}t;
\end{eqnarray*}
here $F_{(Y_{1},\ldots,Y_{n})}$ is the distribution of
$(Y_{1},\ldots,Y_{n})$ which, by definition,  can be recovered by the
product rule $F_{(Y_{1},\ldots,Y_{n})}(y_{1},\ldots, y_{n})=
F_{Y_{1}}(y_{1})F_{Y_{2}|Y_{1}}(y_{2};y_{1})\cdots
F_{Y_{n}|Y_{1},\ldots,Y_{n-1} }(y_{n};y_{1},\ldots,\allowbreak y_{n-1})$
with $F_{1}=A_{0}$ and
\[\label{eqpredpolya}
F_{Y_{j}|Y_{1},\ldots,Y_{j-1} }(y;y_{1},\ldots,y_{j-1})
=\frac{a}{a+j-1}A_{0}(y)+\frac{1}{a+j-1}\sum_{i=1}^{j-1}\mathbbl{1}_{(-\infty,y]}(y_{i}).
\]
However, some remarks on the asymptotic behavior of $\varepsilon(n)$ can be made under
suitable conditions. Since,
\[
\frac{1}{t}f_{T^{(n)}}\biggl(\frac{z-(1-t)x}{t}\biggr)\geq f_{T^{(n)}}\biggl(\frac{z-(1-t)x}{t}\biggr),
\]
if the support of $Y_{i}$'s is bounded (for instance equal to $[0,1]$) and the derivative of $f_{T^{(n)}}$ is bounded by some constant $K$, then, by Taylor expansion of
$f_{T^{(n)}}$, we have
%
\begin{eqnarray}\label{eqqualeps}
\varepsilon(n)&=&\frac{1}{B(a,n)}\int_{0}^{1}(1-t)^{a-1}t^{n}
\biggl(\int_{0}^{1}\frac{1}{t}f_{T^{(n)}}(z/t)\,\mathrm{d}z\biggr)\,\mathrm{d}t\nonumber\\
&&{}-\frac{1}{B(a,n)}\sup_{x\in C^{(n)}}x\int_{0}^{1}
\biggl(\int_{0}^{1}(1-t)^{a}t^{n-2}f^{'}_{T^{(n)}}\biggl(\frac{z_{x,t}}{t}\biggr)\,\mathrm{d}t\biggr)\,\mathrm{d}z\nonumber\\ [-8pt]\\ [-8pt]
&\geq&\frac{n}{a+n}-KK^{(n)}(\lambda)\int_{0}^{1}(1-t)^{a}t^{n-2}\,\mathrm{d}t\nonumber\\
&=&\frac{n}{a+n}-\frac{aKK^{(n)}(\lambda)}{n-1}.\nonumber
\end{eqnarray}
For a large enough $n_{0}$, if we fix $\lambda$ equal to some positive constant $C$
which is grater than $a/(a+n)$ for all $n>n_{0}$, then $K^{(n)}(\lambda)$ is bounded above by
\[
\frac{1-C+\E[Y_{1,1}]}{C-a/(a+n_{0})}.
\]
The second term in \eqref{eqqualeps} is negligible with respect to the first term,
which increase as~$n$ increases. As we mentioned, when the support of the $Y_{i}$'s is
bounded, from Theorem~16.2.4 in \cite{Mey93} 
it follows that the error in total variation between the $m$th transition
probability of the Markov chain $\{M_{m}^{(n)},m\geq0\}$ and the limit distribution $\mathscr{M}$ is
less than $(1-\varepsilon(n))^{m}$. This error decreases for $n$ increasing greater
than $n_0$.

So far we have provided only some qualitative features on the rate of convergence; however,
the derivation of the explicit bound of the rate of convergence of $M_{m}^{(n)}$ to $M$ for
each fixed $n$, via $p_0^{(n)}$ and $\varepsilon(n)$, is still an open problem.
Some examples
confirm our conjecture that the convergence of the Markov chain $\{M_{m}^{(n)},m\geq0\}$
improves as~$n$ increases. Nonetheless, we must point out that
simulating the innovation term $T^{(n)}$ for~$n$ larger than 1 will be more computationally expensive, and also that this cost will be increasing
as~$n$ increases. In fact, if $n$ is greater that one, $2n-1$ more r.v.s
must be drawn  at each iteration  of the Markov chain ($n-1$ more from the  P\'olya sequence
and~$n$ more from the finite-dimensional Dirichlet distribution).
Moreover, we compared the total user times of the R function simulating $\{M_m^{(n)},m=0,\ldots, 500\}$.
We found that all these times were small, of course depending on
$\alpha_0$,  but not on the total mass parameter~$a$ (all the other values being fixed).
The total user times when $n=2$ were about $50\%$ greater than those for $n=1$, while
they were about 5, 10 and 50 times greater when $n=10, 20$ and 100, respectively, for
a number of total iterations equal to 500. From the following
examples, we found that values of $n$ between $2$ and $20$ are a good choice between
a fast rate of convergence and a moderate computational cost.

\begin{exe}\label{ex1}
Let $\alpha_{0}$ be a Uniform distribution on $(0,1)$ and let $a$ be the total mass.
In this case $\E[|Y_{1,j}|]=1/2$ so that for any fixed integer $n$,
the chain will be geometrically ergodic; moreover, it can be proved that $(0,1)$ is small
so that the chain is uniformly ergodic. When $a=1$, \cite{Gug01}
showed
that the convergence of $\{M_{m},m\geq0\}$ is very good and there is no need to consider
the chain
with $n>1$. We consider the cases $a=10$, $50$ and $100$, and for each of them we run
the chain for $n=1$, $2$, $10$ and $20$. We found that the trace plots do not depend on
the initial values. In Figure \ref{figfigure1bis},
we give the trace plots of~$M_{m}^{(n)}$ when
$M_{0}^{(n)}=0$.  Observe that convergence improves as $n$ increases for any fixed value
of $a$; however the improvement is more glaring from the graph for large $a$.
When $a=100$ the convergence of
the chain for $n=1$ seems to occur
at about $m=350$, while for $n=20$ the convergence is at about a value between $50$ and $75$.
For these values of $n$, the total user times to reach convergence was 0.038 seconds for the
former, and 0.066 seconds for the latter. Moreover, the total user times to simulate
500 iterations of $\{M_{m}^{(n)},m\geq0\}$ were 0.05, 0.071, 0.226, 0.429, 2.299 seconds
when $n=1, 2, 10, 20, 100$, respectively.\looseness=1

\begin{figure}

\includegraphics{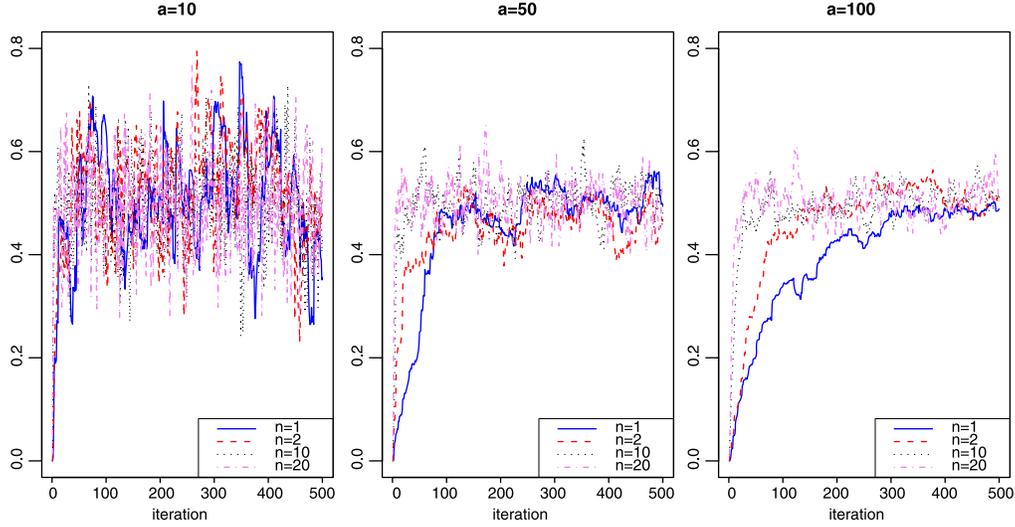}

\caption{Traceplots of the Markov chain $\{M_{m}^{(n)},m\geq0\}$ with $\alpha_{0}$ the Uniform distribution on $(0,1)$, $a=10,50,100$
and $n=1$ (solid blue line), $n=2$ (dashed red line), $n=10$ (dotted black line) and $n=20$ (dot-dashed violet line).}
\label{figfigure1bis}
\end{figure}
\end{exe}

This behaviour is confirmed in the next example, where
the support of the measure $\alpha$ is assumed to be
unbounded.

\begin{figure}

\includegraphics{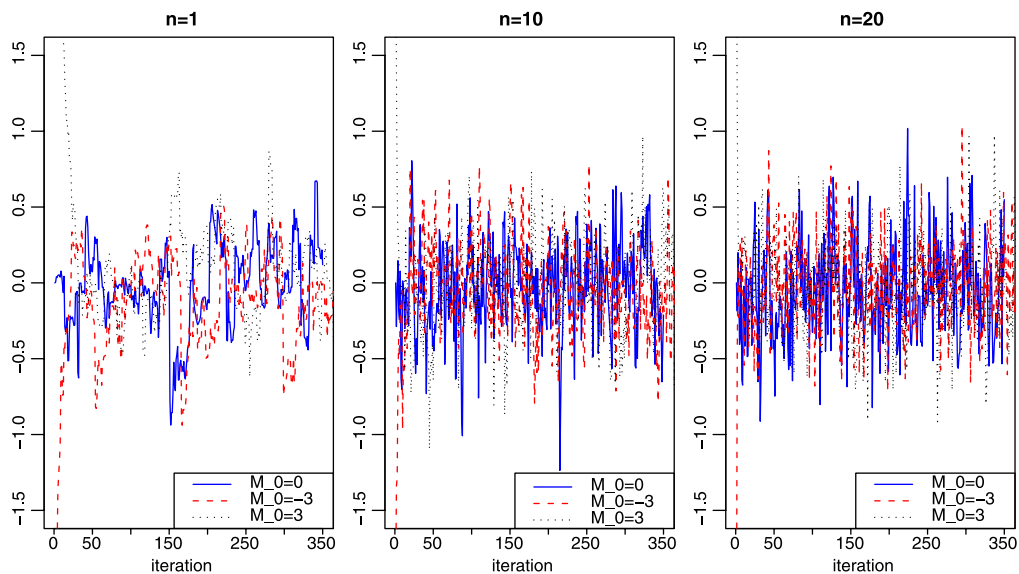}

 \caption{Traceplots of the Markov chain $\{M_{m}^{(n)},m\geq0\}$ with $\alpha_{0}$ the Gaussian distribution with parameter
$(0,1)$, $a=10$, $n=1,10,20$ and $M^{(1)}_{0}=-3$ (dashed red line), $M^{(1)}_{0}=0$ (solid blue line) and $M^{(1)}_{0}=3$ (dotted black line).}
\label{figfigura2}
\end{figure}

\begin{exe}\label{ex2}
Let $\alpha_{0}$ be a Gaussian distribution with parameter $(0,1)$ and let $a=10$.
The behavior of ${M_{m},m\geq0}$ has been considered
in \cite{Gug01}. 
Figure~\ref{figfigura2} displays the trace plots of
${M_{m}^{(n)},m\geq0}$ for three different initial values
($M_{0}^{(n)}=-3,0,3$), with $n=1,10,20$.
Also in this case, it is clear that the convergence improves as $n$ increases.
As far as the total user times are concerned, we drew similar conclusions than in
Example~\ref{ex1}.
\end{exe}

The next is an example in the mixture models context.

\begin{exe}\label{ex3}
Let us consider a Gaussian kernel $k(x,\theta)$ with unknown mean $\theta$ and known
variance equal to 1. If we consider the random density $f(x)=\int_{\mathbb{R}} k(x,\theta)\,\mathrm{d}P(\theta)$,
where $P$ is a Dirichlet process with parameter $\alpha$, then, for any fixed $x$, $f(x)$ is
a random Dirichlet mean. Therefore, if we consider the measure-valued Markov chain
$\{ P_m^{(n)}, m\geq 0\}$ defined recursively as in \eqref{ftchain4}, we define a
sequence of random densities $\{ f_m^{(n)}(x), m\geq 0\}$, where
$f_m^{(n)}(x):=\int_{\mathbb{R}} k(x,\theta)\,\mathrm{d}P_m^{(n)}(\theta)=\theta_m \sum_1^nq_{m,i}^{(n)} k(x,Y_{m,i})+(1-\theta_m)f_{m-1}^{(n)}(x)$.
In each panel of Figure~\ref{figfigura3}, we drew $f_m^{(n)}(x)$
for different values of $m$ when $n$ is fixed. In particular, we fixed $\alpha_0$ to be
a Gaussian distribution with parameter $(0,1)$, and let $a=100$; in this case,
since the ``variance'' of $P$ is small, the mean density
$E[f](x)=\int_{\mathbb{R}} k(x,\theta)\alpha_0(\mathrm{d}\theta)$ (Gaussian with parameter $(0,2)$)
will be very close to the random function
$f(x)$,  so that it can be considered an
approximation of the ``true'' density $f(x)$. From the plots, it is clear
that the convergence improves as $n$ increases: when $n=1$, only $f_{1000}^{(1)}(x)$ is
close enough to the  mean density  $E[f](x)$, while, if $n=20$,   $f_{100}^{(20)}(x)$,
as well as the successive iterations, is a~good approximation of  $f(x)$.
%
\begin{figure}

\includegraphics{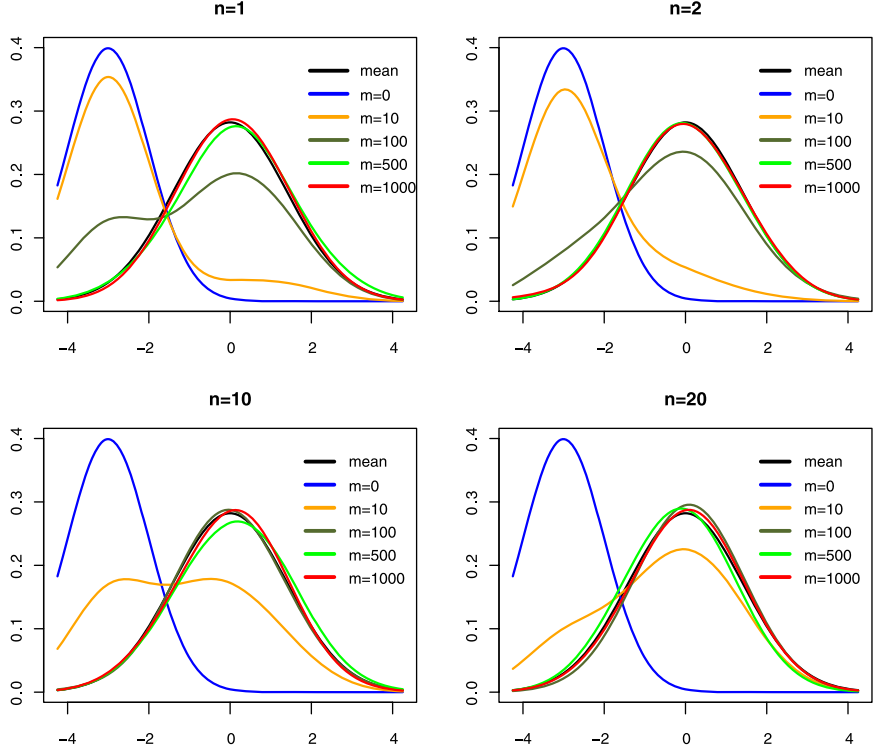}

 \caption{Plots of  $f_{m}^{(n)}$ as in Example~\protect\ref{ex3};  $a=100$,
$\alpha_{0}={\mathcal N}(0,1)$, $n=1,2,10,20$ and $f^{(n)}_0={\mathcal N}(\cdot;-3,1)$.}
\label{figfigura3}
\end{figure}
In any case, observe that even if the ``true'' density~$f(x)$ is unknown, as
when $a=1$, the improvement (as $n$ increases) is clear as well; see
Figure~\ref{figfigura4}, where 5 draws of $f_{m}^{(n)}(x)$, $m=1,100,1000$, are
plotted for different values of $n$.

\begin{figure}

\includegraphics{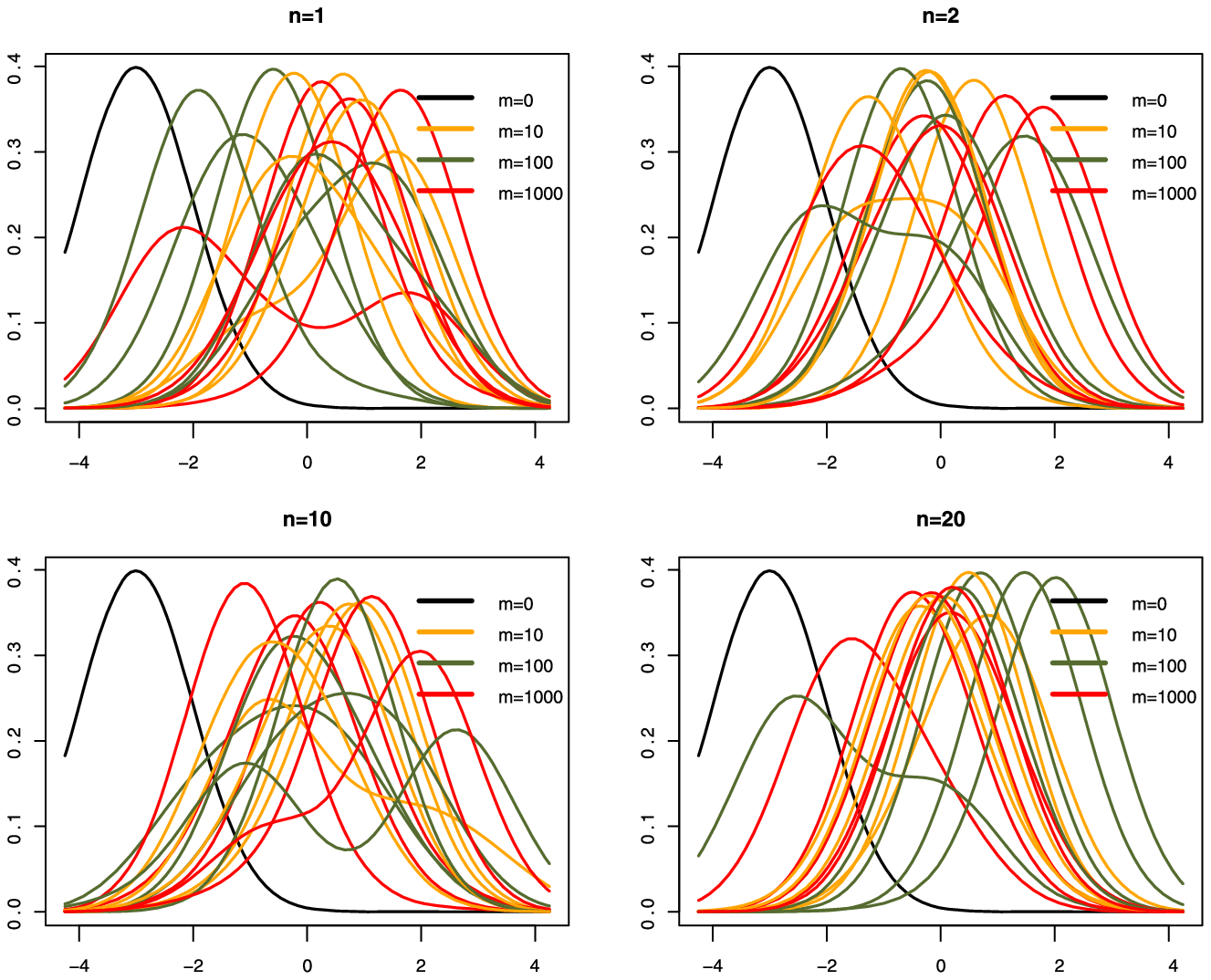}

 \caption{Plots of 5 draws from $f_{m}^{(n)}$ as in Example~\protect\ref{ex3}; $a=1$,
$\alpha_{0}={\mathcal N}(0,1)$, $n=1,2,10,20$ and $f^{(n)}_0={\mathcal N}(\cdot;-3,1)$.}
\label{figfigura4}
\end{figure}
\end{exe}

\section{Concluding remarks and developments}\label{sec4}
The paper \cite{Fei89} 
constitutes, as far as we know, the first work
highlighting the interplay between Bayesian nonparametrics on the one side and the theory of
measure-valued Markov chains on the other. In the present paper, we have further studied
such interplay by introducing and investigating a new measure-valued Markov chain
$\{Q_{n},n\ge1\}$ defined via exchangeable sequences of r.v.s. Two applications
related to Bayesian nonparametrics have been considered: the first gives evidence that $\{Q_{n},n\ge1\}$
is strictly related to the Newton's algorithm of the mixture of Dirichlet process model, while the second shows how
$\{Q_{n},n\geq1\}$ can be applied in order to a define a generalization of the Feigin--Tweedie
Markov chain.

An interesting development consists in investigating whether there are any new
applications related to the Feigin--Tweedie Markov chain apart from the well-known
application in the field of functional linear functionals of the Dirichlet process
(see, e.g., \cite{Gug02}). 
The
proposed generalization of the Feigin--Tweedie Markov chain represents a large class of
measure-valued Markov chains $\{\{P_{m}^{(n)},m\geq 0\},n\in\mathbb{N}\}$ maintaining all the same
properties of the Feigin--Tweedie Markov chain; in other terms, we have increased the flexibility of the
Feigin--Tweedie Markov chain via a further parameter $n\in\mathbb{N}$.
We believe that a number of different interpretations for $n$ can be investigated
in order to extend the applicability of the Feigin--Tweedie Markov chain.

In this respect, an intuitive and simple
extension is related to the problem of defining a (bivariate) vector of measure-valued
Markov chains  $\{ (P_m^{(n_{1})},P_m^{(n_{2})} )  ,m\geq 0\}$,
where, for each fixed $m$, $(P_m^{(n_{1})},P_m^{(n_{2})} ) $ is a vector of dependent
random probabilities, $n_{1}<n_{2}$ being fixed positive integers.
Marginally, the two sequences  $\{P_{m}^{(n_{1})},m\geq0\}$ and $\{P_{m}^{(n_{2})},
m\geq0\}$ are defined via the recursive identity \eqref{ftchain4}; the
dependence is achieved using the same  P\'olya sequence $(Y_{m,1}^{(n_2)},\ldots,$ $Y_{m,n_1}^{(n_2)},
\ldots, Y_{m,n_2}^{(n_2)} )$ and assuming dependence in
$(\theta_m^{(n_2)},\theta_m^{(n_1)})$ or between
$(q_{m,1}^{(n_{2})},\ldots,q_{m,n_2}^{(n_{2})}  )$ and
$(q_{m,1}^{(n_{1})},\ldots,q_{m,n_1}^{(n_{1})})$. For instance, if,
for each $m$, $Z_m, Z_m^{(1)},\ldots,Z_{m}^{(n_{1})},\ldots,Z_m^{(n_2)}$ are
independent r.v.s, $Z_m$ with an Exponential distribution with parameter $a$,
$Z_{m}^{(i)}$  with an Exponential distribution with parameter 1, we could define
$\theta_m^{(n_2)}:={ {\sum_{i=1}^{n_{2}}Z_{m}^{(i)} }/
{( Z_{m}+\sum_{i=1}^{n_{2}}Z_{m}^{(i)})} }$,
$\theta_m^{(n_1)} := {\sum_{i=1}^{n_{1}}Z_{m}^{(i)} }/ ( Z_{m} +\sum_{i=1}^{n_{1}}Z_{m}^{(i)}) $.
Of course, the dependence is related to the difference
between $n_1$ and $n_2$. Work on this is ongoing.

\begin{appendix}
\section*{Appendix}
\subsection{\texorpdfstring{Weak convergence for the Markov chain $\{Q_{n},n\geq1\}$}{Weak convergence for the Markov chain \{Q_{n}, n >= 1\}}}\label{appa1}
\renewcommand{\theequation}{\arabic{equation}}

A proof of the weak convergence of the sequence $Q_n=\sum_1^n q_i^{(n)}\delta_{Y_i}$
of r.p.m.s on $\mathbb{X}$  to a~Dirichlet process $P$ is provided
here, when  the $Y_j$'s are a 
P\'olya sequence with
parameter measure $\alpha$. The result automatically follows from Theorem \ref{asconv}(ii), but this proof
is interesting \textit{per se}, since we use a combinatorial technique;
moreover an explicit expression for the
moment of order $(r_{1},\ldots,r_{k}$) of the $k$-dimensional P\'olya distribution is
obtained.

\begin{prop*}
Let $Q_n$ defined in \eqref{eqPrec}, where
$\{Y_{j},j\geq1\}$ are a 
P\'olya sequence with
parameter $\alpha$. Then
\[
Q_n\Rightarrow P,
\]
where $P$ is a Dirichlet process on $\mathbb{X}$ with parameter $\alpha$.
\end{prop*}

\begin{pf}
By Theorem 4.2 in \cite{Kal83}, 
it is sufficient to prove that for
any measurable partition $B_{1},\ldots,B_{k}$ of $\mathbb{X}$,
\[
(Q_{n}(B_1),\ldots,Q_n(B_k)) \Rightarrow (P(B_1),\ldots,P(B_k)),
\]
characterizing the distribution of the limit. For any given
measurable partition $B_{1},\ldots,B_{k}$ of $\mathbb{X}$, by conditioning on $Y_1,\ldots,Y_n$, it can
be checked that $(Q_{n}(B_{1}),\ldots,Q_{n}(B_{k-1}))$ is distributed according to a
Dirichlet distribution with empirical parameter
$(\sum_{1\leq i\leq n}\delta_{Y_{i}}(B_{1}),\allowbreak \ldots,\sum_{1\leq i\leq n}\delta_{Y_{i}}(B_{k-1}),
\sum_{1\leq i\leq n}\delta_{Y_{i}}(B_{k}))$, and
%
\begin{eqnarray}\label{Polyapartition}
&&\mathbb{P}(\#\{i\dvt Y_i\in B_1\}=j_1,\ldots,\#\{i\dvt Y_i\in B_k\}=j_k
)\nonumber\\ [-8pt]\\ [-8pt]
&&\quad={n\choose j_{1}\cdots j_{k}}\frac{(\alpha(B_{1}))_{j_{1}\uparrow1} \cdots(\alpha(B_{k}))_{j_{k}\uparrow1}}{(a)_{n\uparrow1
}},\nonumber
\end{eqnarray}
where $(j_1,\ldots,j_k)\in\mathcal{D}^{(0)}_{k,n}$ with $\mathcal{D}^{(0)}_{k,n}:=\{(j_{1},\ldots,j_{k})\in\{0,\ldots,n\}^{k}
:\sum_{1\leq i\leq k}j_{i}=n\}$.
For any $k$-uple of nonnegative integers $(r_1,\ldots,r_k)$, we are going to compute
the limit, for $n\rightarrow+\infty$, of the moment
%
\begin{eqnarray}\label{mommisto-gener-Hn}
&&\hspace*{-5pt}\E\Biggl[(Q_n(B_1))^{r_1}\cdots(Q_n(B_{k-1}))^{r_{k-1}}
\Biggl( 1-\sum_{i=1}^{k-1}Q_{n}(B_{i})\Biggr)^{r_k}\Biggr]\nonumber\\ [-8pt]\\ [-8pt]
&&\hspace*{-5pt}\quad=
\sum_{(j_{1},\ldots,j_{k})\in\mathcal{D}^{(0)}_{k,n} }{n\choose j_{1}\cdots j_{k} }
\frac{(\alpha(B_{1}))_{j_{1}\uparrow1} \cdots(\alpha(B_{k}))_{j_{k}\uparrow1}}{(a)_{n\uparrow1 }}
\frac{(j_{1})_{r_{1}\uparrow1}\cdots(j_{k})_{r_{k}\uparrow1} }{(n)_{(r_{1}+\cdots+r_{k})\uparrow1}},\nonumber
\end{eqnarray}
where in general $(x)_{n\uparrow\alpha}$ denotes the Pochhammer symbol for the $n$th factorial
power of~$x$ with increment $\alpha$, that is $(x)_{n\uparrow\alpha}:=\prod_{0\leq i\leq n-1}(x+i\alpha)$.
We will show that, as $n\rightarrow +\infty$,
\begin{eqnarray*}
&&\E\Biggl[(Q_n(B_1))^{r_1}\cdots(Q_n(B_{k-1}))^{r_{k-1}}
\Biggl( 1-\sum_{i=1}^{k-1}Q_{n}(B_{i})\Biggr)^{r_k}\Biggr]\\
&&\quad\rightarrow\E\Biggl[(P(B_1)\big)^{r_1}\cdots(P(B_{k-1}))^{r_{k-1}}
\Biggl( 1-\sum_{i=1}^{k-1}P(B_{i})\Biggr)^{r_k}\Biggr],
\end{eqnarray*}
where $P$ is a Dirichlet process on $\X$ with parameter measure $\alpha$, that is, the r.v.
$(P(B_1),\ldots, P(B_{k-1}))$ has Dirichlet distribution with parameter
$(\alpha(B_1),\ldots,\alpha(B_{k-1}),\allowbreak \alpha(B_k))$. This will be sufficient to characterize the distribution of the limit $Q^*$,
because of the boundedness of the support of the limit distribution.
First of all, we prove the convergence for $k=2$, which corresponds to the one-dimensional case. In particular, we have
\begin{eqnarray*}
&&\E\bigl[(Q_n(B_1))^{r_1}
\bigl( 1-(Q_n(B_1))\bigr)^{r_2}\bigr]\\
&&\quad= \frac{1}{(n)_{(r_{1}+r_{2})\uparrow1}}\sum_{(j_{1},j_{2})\in\mathcal{D}^{(0)}_{2,n}  }
{n\choose j_{1},j_{2} }\frac{(\alpha(B_{1}))_{j_{1}\uparrow1}(\alpha(B_{2}))_{j_{2}\uparrow1}}
{(a)_{n\uparrow1 }}(j_{1})_{r_{1}\uparrow1}(j_{2})_{r_{2}\uparrow1}\\
&&\quad=\frac{1}{(n)_{(r_{1}+r_{2})\uparrow1}}\sum_{t_{1}=1}^{r_{1}}|s(r_{1},t_{1})|\sum_{s_{1}=1}^{t_{1}}S(t_{1},s_{1})
\sum_{t_{2}=1}^{r_{2}}|s(r_{2},t_{2})|\sum_{s_{2}=1}^{t_{2}}S(t_{2},s_{2})\\
&&\qquad{}\times\sum_{(j_{1},j_{2})\in\mathcal{D}^{(0)}_{2,n} }{n\choose j_{1},j_{2} }\frac{(\alpha(B_{1}))_{j_{1}
\uparrow1}(\alpha(B_{2}))_{j_{2}\uparrow1}}{(a)_{n\uparrow1
}}(j_{1})_{s_{1}\downarrow1}(j_{2})_{s_{2}\downarrow1},
\end{eqnarray*}
where $(x)_{n\downarrow 1}:=(-1)^{-n}(-x)_{n\uparrow1}$ and $s(\cdot,\cdot)$ and
$S(\cdot,\cdot)$ are the Stirling number of the first and second kind, respectively.
Let us consider the following numbers, where $s_1$, $s_2$ are nonnegative integers and $n=1,2,\ldots,$
\[
C_{n}^{(s_{1},s_{2})}:=\sum_{(j_{1},j_{2})\in\mathcal{D}^{(0)}_{2,n} }
{n\choose j_{1},j_{2} }\frac{(\alpha(B_{1}))_{j_{1}\uparrow1}(\alpha(B_{2}))_{j_{2}\uparrow1}}
{(a)_{n\uparrow1 }}(j_{1})_{s_{1}\downarrow1}(j_{2})_{s_{2}\downarrow1},
\]
and prove they satisfy a recursive relation. In particular,
\begin{eqnarray*}
C_{n+1}^{(s_{1},s_{2})}& =& \sum_{(j_{1},j_{2})\in\mathcal{D}^{(0)}_{2,n+1} }
{n+1\choose j_{1},j_{2} }\frac{(\alpha(B_{1}))_{j_{1}\uparrow1}(\alpha(B_{2}))_{j_{2}
\uparrow1}}{(a)_{(n+1)\uparrow1 }}(j_{1})_{s_{1}\downarrow1}(j_{2})_{s_{2}\downarrow1}\\
&=&\sum_{j_{1}=0}^{n}{n+1\choose j_{1}+1}\frac{(\alpha(B_{1}))_{(j_{1}+1)\uparrow1}(\alpha(B_{2}))_{(n-j_{1})\uparrow1}}{(a)_{(n+1)\uparrow1 }}(j_{1}+1)_{s_{1}\downarrow1}(n-j_{1})_{s_{2}\downarrow1}\\
&=&\frac{(n+1)(\alpha(B_{1})+s_{1}-1)}{(a+n)}C_{n}^{(s_{1}-1,s_{2})}+\frac{n+1}{(a+n)}C_{n}^{(s_{1},s_{2})},
\end{eqnarray*}
so that the following recursive equation holds
%
\begin{equation}\label{rec1}
C_{n+1}^{(s_{1},s_{2})}=\frac{(n+1)(\alpha(B_{1})+s_{1}-1)}{(a+n)}C_{n}^{(s_{1}-1,s_{2})}+\frac{n+1}{(a+n) }C_{n}^{(s_{1},s_{2})}.
\end{equation}
Therefore, starting from $C_n^{(0,0)}=1$,
$C_n^{(0,1)}={ n\alpha(B_2)/a}$, we have
\begin{eqnarray*}
C_{n}^{(1,1)} &=& \sum_{i=0}^{n-1}\frac{(i+1)\alpha(B_{1})}{(a+i)}C_{i}^{(0,1)}\prod_{j=i+1}^{n-1}
\frac{j+1}{a+j}\\
&=&\frac{\Gamma(n+1)\alpha(B_{1}) }{\Gamma(a+n)}\sum_{i=0}^{n-1}\frac{\Gamma(a+i)}{\Gamma(i+1)}C_{i}^{(0,1)}\\
&=& \frac{\Gamma(n+1)\alpha(B_{1}) }{\Gamma(a+n)}\sum_{i=0}^{n-1}\frac{\Gamma(a+i)}{\Gamma(i+1)}\frac{i\alpha(B_{2})}{a}
=\frac{\alpha(B_{1})\alpha(B_{2})(n)_{2\downarrow1}}{(a)_{2\uparrow1}},
\end{eqnarray*}
and by \eqref{rec1} we obtain $C_{n}^{(s_{1},s_{2})}=(\alpha(B_{1}))_{s_{1}\uparrow1}  (\alpha(B_{2}))_{s_{2}\uparrow1} (n)_{(s_{1}+s_{2})\downarrow1}/(a)_{(s_{1}+s_{2})
\uparrow1}$.
Thus,
\begin{eqnarray*}
&&\lim_{n\rightarrow+\infty} \E\bigl[(H_n(B_1))^{r_1}
\bigl( 1-(H_n(B_1))\bigr)^{r_2}\bigr]\\
&&\quad =\sum_{t_{1}=0}^{r_{1}}|s(r_{1},t_{1})|\sum_{s_{1}=0}^{t_{1}}S(t_{1},s_{1})\sum_{t_{2}=0}^{r_{2}}|s(r_{2},t_{2})|\sum_{s_{2}=0}^{t_{2}}S(t_{2},s_{2})\lim_{n\rightarrow+\infty}\frac{C_{n}^{(s_{1},s_{2})}}{(n)_{(r_{1}+r_{2})\uparrow1}}  \\
&&\quad =|s(r_{1},r_{1})|S(r_{1},r_{1})|s(r_{2},r_{2})|S(r_{2},r_{2})\frac{(\alpha(B_{1}))_{r_{1}\uparrow1}
(\alpha(B_{2}))_{r_{2}\uparrow1}  }{(a)_{(r_{1}+r_{2})\uparrow1}}\\
&&\quad=\frac{(\alpha(B_{1}))_{r_{1}\uparrow1}(\alpha(B_{2}))_{r_{2}\uparrow1}  }{(a)_{(r_{1}+r_{2})\uparrow1}}.
\end{eqnarray*}
The last expression is exactly $\E[ (P(B_1))^{r_1} (1-P(B_1))^{r_2}]$, where
$P(B_1)$ has Beta distribution with parameter $(\alpha(B_1),\alpha(B_2))$.

This proof can be easily generalized to the case $k>2$.
Analogously to the one-dimensional case, we can write
\begin{eqnarray*}
&&\E\Biggl[(Q_n(B_1))^{r_1}\cdots(Q_n(B_{k-1}))^{r_{k-1}}\Biggl( 1-\sum_{i=1}^{k-1}Q_n(B_i)\Biggr)^{r_k}\Biggr]\\
&&\quad=\frac{1}{(n)_{(r_{1}+\cdots+r_{k})\uparrow1}}\sum_{(j_{1},\ldots,j_{k})\in \mathcal{D}^{(0)}_{k,n} }{n\choose j_{1}\cdots j_{k} }\frac{(\alpha(B_{1}))_{j_{1}\uparrow1}\cdots(\alpha(B_{k}))_{j_{k}\uparrow1}}{(a)_{n\uparrow1 }}
(j_{1})_{r_{1}\uparrow1}\cdots(j_{k})_{r_{k}\uparrow1}\\
&&\quad=\frac{1}{(n)_{(r_{1}+\cdots+r_{k})\uparrow1}}\sum_{t_{1}=0}^{r_{1}}|s(r_{1},t_{1})|\sum_{s_{1}=0}^{t_{1}}S(t_{1},s_{1})
\cdots\sum_{t_{k}=0}^{r_{k}}|s(r_{k},t_{k})|\sum_{s_{k}=0}^{t_{k}}S(t_{k},s_{k})\\
&&\qquad{} \times\sum_{(j_{1},\ldots,j_{k})\in\mathcal{D}^{(0)}_{k,n} }
{n\choose j_{1}\cdots j_{k} }\frac{(\alpha(B_{1}))_{j_{1}\uparrow1}\cdots(\alpha(B_{k}))_{j_{k}\uparrow1}}{(a)_{n\uparrow1 }}
(j_{1})_{s_{1}\downarrow1}\cdots(j_{k})_{s_{k}\downarrow1}
\end{eqnarray*}
and, as before, define
\[
C_{n}^{(s_{1},\ldots,s_{k})}:=\sum_{(j_{1},\ldots,j_{k})\in\mathcal{D}^{(0)}_{k,n} }
{n\choose j_{1}\cdots j_{k} }\frac{(\alpha(B_{1}))_{j_{1}\uparrow1}\cdots(\alpha(B_{k}))_{j_{k}\uparrow1}}{(a)_{n\uparrow1 }}
\times(j_{1})_{s_{1}\downarrow1}\cdots(j_{k})_{s_{k}\downarrow1}
\]
and prove they satisfy a recursive relation.
Observe that $C_{n}^{(s_{1},\ldots,s_{k})}/{(n)_{(r_{1}+\cdots+r_{k})\uparrow1}}$ is the moment of order $(r_1,\ldots,r_k)$ of  the $k$-dimensional P\'olya distribution by
definition.
Therefore, for $k>2$
\begin{eqnarray*}
C_{n}^{(s_{1},\ldots,s_{k})}&=& \sum_{j_k=0}^{n}{n\choose j_k}\frac{(\alpha(B_{k}))_{j_k\uparrow1}(a-\alpha(B_{k}))_{(n-j_{k})\uparrow1}
(j_{k})_{s_{k}\downarrow1}}{(a)_{n\uparrow1 }}\\[-2pt]
&&\hphantom{\sum_{j_k=0}^{n}}{}\times\sum_{(j_{1},\ldots,j_{k-1})\in\mathcal{D}^{(0)}_{k-1,n-j_{k}} }{n-j_{k}\choose j_{1}\cdots j_{k-1} }
\frac{(\alpha(B_{1}))_{j_{1}\uparrow1}\cdots(\alpha(B_{k-1}))_{j_{k-1}\uparrow1}}{(a-\alpha(B_{k}))_{(n-j_{k})\uparrow1
}}\\[-2pt]
&&\hspace*{10pt}\hphantom{\sum_{j_k=0}^{n}\sum_{(j_{1},\ldots,j_{k-1})\in\mathcal{D}^{(0)}_{k-1,n-j_{k}} }}{}\times(j_{1})_{s_{1}\downarrow1}\cdots(j_{k-1})_{s_{k-1}\downarrow1}\\[-2pt]
&=& \sum_{j_{k}=0}^{n}{n\choose j_{k}}\frac{(\alpha(B_{k}))_{j_{k}\uparrow1}(a-\alpha(B_{k}))_{(n-j_{k})\uparrow1} (j_{k})_{s_{k}\downarrow1}}
{(a)_{n\uparrow1 }}\\[-2pt]
&&\hphantom{\sum_{j_{k}=0}^{n}}{}\times\frac{(\alpha(B_{1}))_{s_{1}\uparrow1}\cdots (\alpha(B_{k-1}))_{s_{k-1}\uparrow1}
(n-j_{k})_{(s_{1}+\cdots+ s_{k-1})\downarrow1}}{(a-\alpha(B_{k}))_{(s_{1}+\cdots+s_{k-1})\uparrow1}},
\end{eqnarray*}
where the last equality follows by induction hypothesis (we already proved the base case $k=2$ in \eqref{rec1}).
Then, following the same steps of the one-dimensional case, we can recover a recursive
equation for $C_{n}^{(s_{1},\ldots,s_{k})}$,
%
\begin{equation}\label{rec2}
C_{n+1}^{(s_{1},\ldots,s_{k})}=\frac{(n+1)(\alpha(B_{k})+s_{k}-1)}{a+n}C_{n}^{(s_{1},
\ldots,s_{k-1},s_{k}-1)}+\frac{n+1}{a+n}C_{n}^{(s_{1},\ldots,s_{k})}.
\end{equation}
Starting from $C_n^{(0,\ldots,0)}=1$, $C_n^{(0,\ldots,0,1,0,\ldots,0)}={n\alpha(B_j)/a} $ and
\begin{eqnarray*}
C_{n}^{(1,\ldots,1)}&=&\sum_{i=0}^{n-1}\frac{(i+1)\alpha(B_{k})}{(a+i)}C_{i}^{(1,\ldots,1,0) }\prod_{j=i+1}^{n-1} \frac{j+1}{a+j} \\[-2pt]
&=&\frac{\Gamma(n+1)}{\Gamma(a+n)}\frac{(\alpha(B_{1}))_{s_{1}\uparrow1}\cdots (\alpha(B_{k-1}))_{s_{p-1}\uparrow1}\alpha(B_{k})}
{(a-\alpha(B_{k}))_{(s_{1}+\cdots+s_{k-1})\uparrow1}}\sum_{i=0}^{n-1}\frac{\Gamma(a+i)}{\Gamma(1+i)}C_{i}^{(1,\ldots,1,0) }\\[-2pt]
&=&\frac{\Gamma(n+1)}{\Gamma(a+n)}\frac{(\alpha(B_{1}))_{s_{1}\uparrow1}\cdots (\alpha(B_{k-1}))_{s_{k-1}\uparrow1}\alpha(B_{k})}
{(a-\alpha(B_{k}))_{(s_{1}+\cdots+s_{k-1})\uparrow1}}\\[-2pt]
&&{} \times(-1)^{k+1}\sum_{i=0}^{n-1}\frac{\Gamma(a+i)}{\Gamma(1+i)}
\frac{(a-\alpha(B_{k}))_{i\uparrow1} (-i)_{(k-1)\uparrow1}}{(a)_{i\uparrow1}}\\[-2pt]
&&\hphantom{{} \times(-1)^{k+1}\sum_{i=0}^{n-1}}{}\times\frac{\Gamma(-a+\alpha(B_{k})-i+1)\Gamma(-k-a+2)}
{\Gamma(-a-i+1)\Gamma(-a+\alpha(B_{k})+2-k)} \\[-2pt]
&=&\frac{\alpha(B_{1})\cdots\alpha(B_{k})(n)_{k\downarrow1}}{(a)_{(k\uparrow1)}},
\end{eqnarray*}
by repeated application of \eqref{rec2}, we obtain
\[
C_{n}^{(s_{1},\ldots,s_{k})}=\frac{(\alpha(B_{1}))_{s_{1}\uparrow1}
\cdots(\alpha(B_{k}))_{s_{k}\uparrow1} (n)_{(s_{1}+\cdots +s_{k})\downarrow1}}{(a)_{(s_{1}+\cdots+s_{k})\uparrow1}}.\vadjust{\goodbreak}
\]
Thus,
\begin{eqnarray*}
&&\lim_{n\rightarrow+\infty}\E\Biggl[(Q_n(B_1))^{r_1}\cdots(Q_n(B_{k-1}))^{r_{k-1}}\Biggl( 1-\sum_{i=1}^{k-1}Q_n(B_i)\Biggr)^{r_p}\Biggr]\\
&&\quad= \sum_{t_{1}=0}^{r_{1}}|s(r_{1},t_{1})|\sum_{s_{1}=0}^{t_{1}}S(t_{1},s_{1})\cdots\sum_{t_{k}=0}^{r_{k}}|s(r_{k},t_{k})|\sum_{s_{k}=0}^{t_{k}}S(t_{k},s_{k})
\lim_{n\rightarrow+\infty}\frac{C_{n}^{(s_{1},\ldots,s_{k})}}{(n)_{(r_{1}+\cdots+r_{k})\uparrow1}}\\
&&\quad=|s(r_{1},r_{1})|S(r_{1},r_{1})\cdots|s(r_{k},r_{k})|S(r_{k},r_{k})
\frac{(\alpha(B_{1}))_{r_{1}\uparrow1}\cdots(\alpha(B_{k}))_{r_{k}\uparrow1}  }{(a)_{(r_{1}+\cdots+r_{k})\uparrow1}}\\
&&\quad=\frac{(\alpha(B_{1}))_{r_{1}\uparrow1}\cdots(\alpha(B_{k}))_{r_{k}\uparrow1}
}{(a)_{(r_{1}+\cdots+r_{k})\uparrow1}}\\
&&\quad=\E\Biggl[(P(B_1))^{r_1}\cdots(P(B_{k-1}))^{r_{k-1}}
\Biggl( 1-\sum_{i=1}^{k-1}P(B_{i})\Biggr)^{r_k}\Biggr],
\end{eqnarray*}
where $P$ is a Dirichlet process with parameter $\alpha$.
\end{pf}

\subsection{Solution of the distributional equation}\label{appa2}
Here, we provide an alternative proof for the solution of the
distributional equation \eqref{ftchain3} introduced by \cite{Fav07}. 

\begin{theorem*}
For any fixed integer $n\geq1$, the distributional equation
\[
P^{(n)}\stackrel{\mathrm{d}}{=}\theta\sum_{i=1}^{n}q_{i}^{(n)}\delta_{Y_{i}}+(1-\theta)P^{(n)}
\]
has the Dirichlet process
with parameter $\alpha$ as its unique solution, assuming the independence between $P^{(n)}$,
$\theta,(q_{1}^{(n)},\ldots,q_{n}^{(n)})$ and $(Y_{1},\ldots,Y_{n})$ in the right-hand side.
\end{theorem*}

\begin{pf}
From Skorohod's theorem, it follows that there exist $n$ independent r.v.s $\xi_{1},\ldots,\xi_{n}$ such that
$\xi_{i}$ has Beta distribution with parameter $(1,n-i)$ for $i=1,\ldots,n$ and $q_{1}^{(n)}=\xi_{1}$ and
$q_{i}^{(n)}=\xi_{i}\prod_{1\leq j\leq i-1 }(1-\xi_{j})$ for $i=2,\ldots,n$; in particular, by a simple transformation of r.v.s,
it follows
that $(q_{1}^{(n)},\ldots,q_{n}^{(n)})$ is distributed according to a~Dirichlet distribution function with
parameter $(1,\ldots,1)$.
Further, since $\xi_{n}=1$ a.s., then $\sum_{1\leq i\leq n}q_{i}^{(n)}=1$ a.s. and it can be
verified by induction that
\[
1-\sum_{i=1}^{j}q_{i}^{(n)}=\prod_{i=1}^{j}(1-\xi_{i}),\qquad j=1,\ldots,n-1.
\]
Let $B_{1},\ldots,B_{k}$ be a measurable partition of $\mathbb{X}$. We first prove that conditionally on $Y_{1},\ldots,Y_{n}$, the
finite dimensional distribution of the r.p.m. $\sum_{1\leq i\leq n}q_{i}^{(n)}\delta_{Y_{i}}$ is the Dirichlet distribution with the empirical parameter
$(\sum_{1\leq i\leq n}\delta_{Y_{i}}(B_{1}),\ldots,\sum_{1\leq i\leq n}\delta_{Y_{i}}(B_{k}))$. Actually, since
\begin{eqnarray*}
&&\Biggl(\Biggl(\sum_{i=1}^{n}q_{i}^{(n)}\delta_{Y_{i}}\Biggr)(\cdot,B_{1}),\ldots,\Biggl(\sum_{i=1}^{n}q_{i}^{(n)}\delta_{Y_{i}}\Biggr)(\cdot,B_{k})\Biggr)\\
&&\quad=\Biggl(\sum_{i=1}^{n}q_{i}^{(n)}\delta_{Y_{i}}(B_{1}),\ldots,\sum_{i=1}^{n}q_{i}^{(n)}\delta_{Y_{i}}(B_{k})\Biggr)\\
&&\quad=\Biggl(\sum_{i:Y_{i}\in B_{1}}q_{i}^{(n)},\ldots,\sum_{i:Y_{i}\in B_{k}}q_{i}^{(n)}\Biggr),
\end{eqnarray*}
then, conditionally on the r.v.s $Y_{1},\ldots,Y_{n}$, the r.v.
$(\sum_{i:Y_{i}\in B_{1}}q_{i}^{(n)}$,
$\ldots,$$\sum_{i:Y_{i}\in B_{k}}q_{i}^{(n)})$ is distributed according to
a Dirichlet distribution with parameter $(n_{1},\ldots,$ $n_{k})$, where $n_{i}=\sum_{1\leq j\leq n}\delta_{Y_{j}}(B_{i})$ for $i=1,\ldots,k$.
Conditionally on $Y_{1},\ldots,Y_{n}$, the finite dimensional distributions of the right-hand side of
\eqref{ftchain3} are Dirichlet with updated parameter
$((\alpha(B_{1})+\sum_{1\leq i\leq n}\delta_{Y_{i}}(B_{1}),\ldots,\alpha(B_{k})+\sum_{1\leq i\leq n}\delta_{Y_{i}}(B_{k}))$.
This argument verifies that the Dirichlet process with parameter $\alpha$ satisfies the distributional equation (\ref{ftchain3}). This solution
is unique by Lemma 3.3 in  \cite{Set94} 
(see also \cite{Ver79}, 
Section~1).
\end{pf}

\subsection{\texorpdfstring{Proofs of the theorems in Section~\protect\ref{secthree}}{Proofs of the theorems in Section 3}}\label{appa3}
\vspace*{-10pt}
\begin{pf*}{Proof of Theorem \ref{thminvariantPn}}
The proof is based on the
``standard'' result that properties (e.g., weak convergence)
of sequences of r.p.m.s can be proved via analogous properties of the sequences of their linear functionals.

First of all, we prove that if $g\dvtx\mathbb{X}\mapsto\mathbb{R}$ is a bounded and
continuous function, then $\{G^{(n)}_{m},m\geq0\}$ with
$G_{m}^{(n)}:=\int_{\mathbb{X}}g(x)P^{(n)}_{m}(\mathrm{d}x)$ is a Markov chain on $\mathbb{R}$ with unique invariant measure $\Pi_{g}$.
From \eqref{ftchain4}, it follows that $\{G^{(n)}_{m},m\geq1\}$ is a Markov chain on $\mathbb{R}$ restricted to
the compact set $[-\sup_{\mathbb{X}}|g(x)| ,\sup_{\mathbb{X}}|g(x)| ]$ and it has
at least one finite invariant measure if it is a weak Feller Markov chain.
In fact, for a fixed $y\in\mathbb{R}$
\begin{eqnarray*}
&&\operatorname{lim\,inf}\limits_{x\rightarrow x^{\ast}}\mathbb{P}\bigl(G^{(n)}_{m}\leq
y\big|G^{(n)}_{m-1}=x\bigr)\\
&&\quad=\operatorname{lim\,inf}\limits_{x\rightarrow x^{\ast}}\mathbb{P}\Biggl(\theta_{m}\sum_{i=1}^{n}q_{m,i}^{(n)}g(Y_{m,i})\leq y-x(1-\theta_{m})\Biggr)\\
&&\quad\ge \int_{(0,1)}\liminf_{x\rightarrow x^{\ast}}\mathbb{P}\Biggl(\sum_{i=1}^{n}q_{m,i}^{(n)}g(Y_{m,i})\leq\frac{y-x(1-z)}{z}\Biggr)
\mathbb{P}(\theta_{m}\in \mathrm{d}z)\\
&&\quad=\int_{(0,1)}\mathbb{P}\Biggl(\sum_{i=1}^{n}q_{m,i}^{(n)}g(Y_{m,i})\leq\frac{y-x^{\ast}(1-z)}{z}\Biggr)\mathbb{P}(\theta_{m}\in
\mathrm{d}z)\\
&&\quad=\mathbb{P}\bigl(G^{(n)}_{m}\leq y\big|G^{(n)}_{m-1}=x^{\ast}\bigr),
\end{eqnarray*}
since the distribution of $\sum_{1\leq i\leq n}q_{m,i}^{(n)}g(Y_{m,i})$ has at most a countable numbers of atoms and $\theta_{m}$ is
absolutely continuous.
From Proposition 4.3 in \cite{Twe76}, 
if we show that
$\{G^{(n)}_{m},\break m\geq0\}$ is $\phi$-irreducible for a finite measure $\phi$, then the Markov chain is
positive
recurrent and the invariant measure $\Pi_g$ is unique. Let us consider the following event $E:=\{Y_{1,1}=Y_{1,2}=\cdots=Y_{1,n}\}$.
Then for a finite measure $\phi$ we have to prove that if $\phi(A)>0$, then $\mathbb{P}(G^{(n)}_{1}\in
A|G^{(n)}_{0})>0$ for any $G^{(n)}_{0}$. We observe that
\begin{eqnarray*}
\mathbb{P}\bigl(G^{(n)}_{1}\in A\big|G^{(n)}_{0}\bigr)
&=&\mathbb{P}\bigl(G^{(n)}_{1}\in A\big|G^{(n)}_{0},E\bigr)\mathbb{P}\bigl(E|G^{(n)}_{0}\bigr)+\mathbb{P}\bigl(G^{(n)}_{1}\in
A\big|G^{(n)}_{0},E^{c}\bigr)\mathbb{P}\bigl(E^{c}|G^{(n)}_{0}\bigr)\\
&\geq& \mathbb{P}\bigl(G^{(n)}_{1}\in A\big|G^{(n)}_{0},E\bigr)\mathbb{P}(E).
\end{eqnarray*}
Therefore, since $\mathbb{P}(E)>0$, using the same argument in Lemma 2 in
\cite{Fei89}, 
we conclude
that $\mathbb{P}(G^{(n)}_{1}\in A|G^{(n)}_{0},E)>0$ for a suitable
measure $\phi$ such that $\phi(A)>0$. Finally, we prove the aperiodicity of $\{G_{m}^{(n)},m\geq0\}$ by contradiction. If the chain is
periodic with period $d>1$, the exist $d$ disjoint sets
$D_{1},\ldots,D_{d}$ such that $\mathbb{P}(G^{(n)}_{m}\in D_{i+1}|G^{(n)}_{m-1}=x)=1$ for all $x\in D_{i}$ and for $i=1,\ldots, d-1$. This implies\vspace*{1pt}
$\mathbb{P}(z\sum_{i=1}^{n}q_{m,i}^{(n)}g(Y_{m,i})+(1-z)x\in D_{i+1})=1$\vspace*{1pt} for almost every $z$ with respect to the Lebesgue measure restricted to $(0,1)$. Thus, $\mathbb{P}(\sum_{i=1}^{n}q_{m,i}^{(n)}g(Y_{m,i})\in D_{i+1})=1$ for  $i=0,\ldots,d-1$.
For generic $\alpha$ and $g$, this is in contradiction with the assumption $d>1$. By Theorem 13.3.4(ii) in \cite{Mey93}, 
$G_{m}^{(n)}$ converges in distribution for $\Pi_g$-almost all starting points
$G_0^{(n)}$.
 In particular, $\{G^{(n)}_{m},m\geq0\}$ converges weakly for $\Pi_{g}$-almost all starting points $G^{(n)}_{0}$.

From the arguments above, it follows that, for all $g$ bounded and continuous,
there exists a r.v. $G$ such that
$G^{(n)}_{m}\Rightarrow G$ as $m\rightarrow+\infty$ for $\Pi_{g}$-almost all starting
points $G^{(n)}_{0}$.
Therefore, for Lemma~5.1 in \cite{Kal83} 
there exists a r.p.m. $P$ such
that $P^{(n)}_{m}\Rightarrow P$ as $m\rightarrow+\infty$ and $G\stackrel{\mathrm{d} }=\int_{\mathbb{X}}g(x)P(\cdot,\mathrm{d}x)$ for all $g\in C(\mathbb{R})$.
This implies that the law of $P$ is an invariant measure for the Markov chain
$\{P^{(n)}_{m},m\geq0\}$. Then, as $m\rightarrow+\infty$,
\[
\int_{\mathbb{X}}g(x)\,\mathrm{d}P^{(n)}_{m}(\cdot,\mathrm{d}x)\Rightarrow\int_{\mathbb{X}}g(x)P(\cdot,\mathrm{d}x)
\]
and the limit is unique for any $g\in C(\mathbb{R})$. Since for any random measure $\zeta_{1}$ and $\zeta_{2}$ we know that
$\zeta_{1}\stackrel{\mathrm{d}}=\zeta_{2}$ if and only if $\int_{\mathbb{X}}g(x)\zeta_{1}(\cdot,\mathrm{d}x)
=\int_{\mathbb{X}}g(x)\zeta_{2}(\cdot,\mathrm{d}x)$ for any
$g\in C(\mathbb{R})$ (see Theorem 3.1. in \cite{Kal83}),
the invariant measure for the Markov chain $\{P^{(n)}_{m},m\geq0\}$ is
unique. By the definition of $\{P^{(n)}_{m},m\geq0\}$, it is straightforward to show that the limit $P$ must satisfy \eqref{ftchain3} so
that $P$ is the Dirichlet process with parameter $\alpha$.
\end{pf*}

\begin{pf*}{Proof of Theorem \ref{thmharrisfunct}}
The proof is a straightforward adaptation of the proof of Theorem~2 in
\cite{Fei89}, 
using
\begin{eqnarray*}
f(G^{(n)}_{1})&=&\log\Biggl(1+\Biggl|\theta_{1}\sum_{i=1}^{n}q_{1,i}^{(n)}g(Y_{1,i})+(1-\theta_{1})G^{(n)}_{0}\Biggr|\Biggr)\\[-2pt]
&\leq& \sum_{i=1}^{n}\log\bigl(1+|g(Y_{1,i})|\bigr)+
\log\bigl(1+(1-\theta_{1})\big|G^{(n)}_{0}|\bigr)
\end{eqnarray*}
instead of their inequality (8).
\end{pf*}

\begin{pf*}{Proof of Theorem \ref{thmgeoergfunct}}
As regards (i), given the definition of stochastically monotone Markov chain, we have that for $z_{1}<z_{2}$, $s\in\mathbb{R}$,
\begin{eqnarray*}
p^{(n)}_{1}(z_{1},(-\infty,s))&=&\mathbb{P}\Biggl(\theta_{1}\sum_{i=1}^{n}q_{1,i}^{(n)}Y_{1,i}+(1-\theta_{1})z_{1}<s\Biggr)\\[-2pt]
&\geq& \mathbb{P}\Biggl(\theta_{1}\sum_{i=1}^{n}q_{1,i}^{(n)}Y_{1,i}+(1-\theta_{1})z_{2}<a\Biggr)\\[-2pt]
&=&p^{(n)}_{1}(z_{2},(-\infty,s)).
\end{eqnarray*}
As far as (ii) is concerned, we first prove that, under condition \eqref{eqfinitmeanalpha0}, the Markov chain
$\{M^{(n)}_{m},m\geq0\}$ satisfies the Foster--Lyapunov condition for the
function $V(x)=1+|x|$. This property implies the geometric ergodicity of the
$\{M^{(n)}_{m},m\geq0\}$ (see \cite{Mey93}, 
Chapter~15). We have
\begin{eqnarray*}
pV(x)&=&\int_{\mathbb{X}} (1+|y|)p(x,\mathrm{d}y)
=1+\E\Biggl[\Biggl|\theta_{1}\sum_{i=1}^{n}q_{1,i}^{(n)}Y_{1,i}+(1-\theta_{1})x\Biggr|\Biggr]\\[-2pt]
&\leq& 1+\E[\theta_{1}]\sum_{i=1}^{n}E\bigl[\bigl|q_{1,i}^{(n)}Y_{1,i}\bigr|\bigr]+|x|\E[1-\theta_{1}]\\[-2pt]
&\leq& 1+\frac{n}{n+a}\sum_{i=1}^{n}\E\bigl[\bigl|q_{1,i}^{(n)}Y_{1,i}\bigr|\bigr]+\frac{a}{n+a}|x|=1+\frac{n}{n+a}\E[|Y_{1,1}|]+\frac{a}{n+a}|x|.
\end{eqnarray*}
Therefore, we are looking for the small set $C^{(n)}$ such that the Foster--Lyapunov
condition holds, that is, a small set
$C^{(n)}$ such that
%
\begin{equation}\label{ftchain15}
1+\frac{n}{n+a}\E[|Y_{1,1}|]+\frac{a}{n+a}|x|\leq\lambda(1+|x|)+b\mathbbl{1}_{C^{(n)}}(x)
\end{equation}
for some constant $b<+\infty$ and $0<\lambda<1$.
If
$C^{(n)}=[-K^{(n)}(\lambda),K^{(n)}(\lambda)]$,
where
\[\label{ftchain17}
K^{(n)}(\lambda):=\frac{1-\lambda+n/(n+a)\E[|Y_{1,1}|]}{\lambda-a/(n+a)},\vadjust{\goodbreak}
\]
then, condition \eqref{ftchain15} holds for all
\[
\lambda\in\biggl(\frac{a}{n+a},1\biggr),\qquad b\geq1-\lambda+\frac{n}{n+a}\E[|Y_{1,1}|].
\]
As in the proof of Theorem~\ref{thminvariantPn}, we can prove that the Markov chain $\{M^{(n)}_{m},m\geq0\}$ is weak Feller; then, since $C^{(n)}$ is a compact set,
it is a small set (see \cite{Twe75}). 
As regards (iii), the proof follows by standard arguments. See, for instance, the proof of Theorem 1 in \cite{Gug01}.
\end{pf*}

\begin{pf*}{Proof of Remark \ref{remjarner}}
As we have already mentioned, the geometric ergodicity follows if a Foster--Lyapunov condition holds.
Let
$V(x)=1+|x|^{s}$; then, if $\E[(1+|\sum_{1\leq i\leq n}q_{1,i}^{(n)}Y_{1,i}|)^{s}]<+\infty$,
it is straightforward to prove that the Foster--Lyapunov condition
$PV(x)\leq\lambda V(x)+b\mathbbl{1}_{\tilde{C}^{(n)}}(x)$ holds for some constant
$b<+\infty$, and $\lambda$ such that
\[
\E[(1-\theta_{1})^{s}]=\frac{\Gamma(a+s)\Gamma(a+n)}{\Gamma(a)\Gamma(a+s+n)}<\lambda<1
\]
and for some compact set $\tilde{C}^{(n)}$. Of course \eqref{eqjar} implies $\E[(1+|\sum_{1\leq i\leq n}q_{1,i}^{(n)}Y_{1,i}|)^{s}]<+\infty;$
in fact, conditioning on
the random number $N$ of distinct values $Y_{1,1}^{\ast},\ldots, Y_{1,N}^{\ast}$ in
$Y_{1,1},\ldots,Y_{1,n}$, $1\leq N\leq n$, we have
\[
\Biggl|\sum_{i=1}^{n}q_{1,i}^{(n)}Y_{1,i}\Biggr|\leq\sum_{i=1}^{N}\tilde{q}_{1,i}^{(n)}|Y^{\ast}_{1,i}|
\leq\max\{|Y_{1,1}^{\ast}|,\ldots,|Y_{1,N}^{\ast}|\}.
\]
Since $\{Y_{1,1}^{\ast},\ldots,Y_{1,N}^{\ast}\}$ are independent and identically distributed according to $\alpha_{0}$,
then
\begin{eqnarray*}
\E\Biggl[\Biggl|\sum_{i=1}^{n}q_{1,i}^{(n)}Y_{1,i}\Biggr|^{s}\Biggr]\leq\int_{0}^{+\infty}y^{s}N(A_{0}(y))^{N-1}\alpha_{0}(\mathrm{d}y)
\leq N\int_{0}^{+\infty}y^{s}\alpha_{0}(\mathrm{d}y)\leq n\E[|Y_{1,i}|^{s}]<+\infty,
\end{eqnarray*}
where $A_{0}$ is the distribution corresponding to the probability measure $\alpha_{0}$,
and this is equivalent to $\E[(1+|\sum_{1\leq i\leq n}q_{1,i}^{(n)}Y_{1,i}|)^{s}]<+\infty$.
\end{pf*}
\end{appendix}

\section*{Acknowledgements}
The authors are very grateful to Patrizia Berti and Pietro Rigo who suggested the proof of Theorem \ref{rigo}, and to
Eugenio Regazzini for helpful discussions. The authors are also grateful to an Associate Editor and a referee for comments
that helped to improve the presentation.
The second author was partially supported by MiUR Grant 2006/134525.


\printhistory

\end{document}